\documentclass[a4paper,12pt]{article}
\usepackage[T2A]{fontenc}
\usepackage[cp1251]{inputenc}
\usepackage[cp1251]{inputenc}
\usepackage[dvips]{graphicx}
\usepackage{amsfonts,amsmath,amssymb,amscd,array, euscript}

\begin{document}

\title{Non-measurable automorphisms of Lie groups
relative to the real- and non-archimedean-valued measures}
\author{S.V. Ludkovsky}
\date{15.10.2007}
\maketitle

\begin{abstract}
In this work the problem about an existence of non-measurable
automorphisms of Lie groups finite and as well infinite dimensional
over the field of real numbers and also over the non-archimedean
local fields is investigated. Non-measurability of automorphisms is
considered relative to real-valued measures and also measures with
values in non-archimedean local fields. Their existence is proved
and a procedure for their construction is given. Their application
for a construction of non-measurable irreducible unitary
representations is demonstrated.
\end{abstract}

\section{Introduction}
Continuous automorphisms of Lie groups were widely studied
\cite{bourgralg,pontr, hewross,fell}, but discontinuous and
non-measurable automorphisms are known substantially less
\cite{bicht,moore,luumn92,lumsb95,lumz2000}. In this work an
explicit direct procedure of a construction of non-measurable
automorphisms of locally compact Lie groups is given and it is
proved that their existence is a local property. Moreover,
non-measurable automorphisms of locally compact Lie algebras are
constructed. Their application for the construction of weakly
non-measurable irreducible unitary representations of locally
compact groups is given. In this article basic necessary facts are
reminded. Besides this non-measurable automorphisms of infinite
dimensional over the real field, as well as non-archimedean fields
of zero characteristic, Lie groups, which are not locally compact
are studied. Non-measurable automorphisms on definite more general
topological groups are investigated. The basic results of the second
section of the paper are Theorems 13, 15, 16, 19, 20 and Corollaries
14, 17, while in Corollary 14, Theorem 19 and \S 13 the specific
features of Lie groups are taken into account.
\par Besides real-valued measures here in the third section
apart from the previous works non-measurability of automorphisms of
totally disconnected topological groups is investigated also for
measures with values in infinite locally compact fields of zero
characteristic with non-archimedean non-trivial multiplicative
norms, that is with values in local fields.

\section{Non-measurable automorphisms of groups relative to
real-valued measures}

\par {\bf 1. Definitions.} For groups $G$ and $S$ a mapping $f: G\to S$
is called a homomorphism, if it preserves the multiplication
operation, that is $f(ab)=f(a)f(b)$ for each $a, b\in G$. If a
homomorphism $f$ is bijective and surjective from $G$ onto $S$,
$f(G)=S$, then the homomorphism $f$ is called the (algebraic)
isomorphism. In the case $G=S$ an isomorphism $f$ is called the
automorphism. \par A homomorphism $f: G\to U(H)$ is called a unitary
representation of the group $G$, if $H$ is the Hilbert or the
unitary space over $\bf C$, $U(H)$ is the unitary group of $H$. In
the particular case $H=\bf C$, that is $U({\bf C})=S^1 = \{ z\in
{\bf C}: |z|=1 \} $, a homomorphism $f: G\to S^1$ is called a
character of the group $G$.
\par If $\sf g$ is a Lie algebra over the field $\bf F$,
then a bijective surjective mapping $\phi : {\sf g}\to \sf g$ we
call an automorphism, if it preserves addition and multiplication:
$\phi (a+b)=\phi (a)+\phi (b)$, $\phi ([a,b])=[\phi (a), \phi (b)]$
for each $a, b\in \sf g$.
\par A topological space $X$ is called compact, if
from each its open covering there is possible to extract a finite
subcovering. A topological space $X$ is called locally compact, if
each its point $x\in X$ has a neighborhood $U$, the closure of which
$\bar U$ is compact (see \cite{eng}; in \cite{pontr} the old
topological terminology is slightly different from the new one
\cite{eng}).
\par For a locally compact Hausdorff topological group
$G$ a non-negative non-trivial $\sigma $-additive measure $\mu $ on
the $\sigma $-algebra ${\cal B}(G)$ of all Borel subsets in $G$ is
called a left (or right) Haar measure, if $\mu (gA) = \mu (A)$ (or
$\mu (Ag)=\mu (A)$ respectively) for each $A\in {\cal B}(G)$ and
$g\in G$.
\par Let a $\sigma $-algebra ${\cal A}(G)= {\cal A}_{\mu }(G)$
be a completion of the Borel $\sigma $-algebra ${\cal B}(G)$ with
the help of subsets $P$ in $G$ such that $P\subset F\in {\cal B}(G)$
and $\mu (F)=0$, where $\mu $ is a non-negative non-trivial $\sigma
$-additive measure on the $\sigma $-algebra ${\cal B}(G)$. An
automorphism $f$ of a topological group $G$ is called $\mu
$-measurable, if $f^{-1}(U)\in {\cal A}(G)$ for each $U\in {\cal
B}(G)$. In the contrary case an automorphism is called $\mu
$-non-measurable.
\par A set $X$ with a $\sigma $-algebra of its subsets $\cal U$
is called a measurable space and it is denoted by $(X, {\cal U})$.
\par A measure $\nu $ is called absolutely continuous  relative to
the measure $\mu $ on the measurable space $(X, {\cal U})$, if from
$\mu (A)=0$ it follows $\nu (A)=0$, where $A\in {\cal U}$. Measures
$\mu $ and $\nu $ are called equivalent, if they are absolutely
continuous relative to each other.
\par {\bf 2. Lemma.} {\it Let $G$ be a separable locally
compact non-compact group, and let $\mu $ be a non-negative Haar
measure on $G$. Then on the one-point compactification $\alpha G$
(considered as the topological space) there exists a finite measure
$\nu $ equivalent to the measure $\mu $.}
\par {\bf Proof.} For each locally compact non-compact topological
space their exists its one-point (Alexandroff) compactification due
to Theorem 3.5.11 \cite{eng}. Take an open neighborhood $U$ of the
unit element in $G$ such that $\mu (U)<\infty $. Since $\mu $ is
non-trivial, then $\mu (U)>0$. Due to the separability of the group
$G$ there exists a countable family $ \{ g_j: j\in {\bf N} \} $ of
elements in $G$ such that $\bigcup_j g_jU=G$, where $gU= \{ z: z=gf,
f\in U \} $. Put
\par $\nu (A) := \sum_{j=1}^{\infty } \mu ((g_jU)\cap A)/2^j$ (1) \\
for each $A\in {\cal A}(G)$ and $\nu ( \{ \alpha \} )=0$, where $ \{
\alpha \} = \alpha G\setminus G$ is the compactification
enlargement. Then the measure $\nu $ is defined on ${\cal A}(\alpha
G)$ and $0<\nu (\alpha G)<\infty $. In view of Formula (1) $\nu
(A)=0$ if and only if $\mu (A)=0$. Therefore, measures $\mu $ and
$\nu $ are equivalent.

\par {\bf 3. Lemma.} {\it If $\mu $ and $\nu $ are two equivalent
$\sigma $-additive measures on ${\cal B}(G)$, where $G$ is a
topological group, then an automorphism $f$ is $\mu $-measurable if
and only if it is $\nu $-measurable.}
\par {\bf Proof.} Since measures $\mu $ and $\nu $
are equivalent and on the same $\sigma $-algebra ${\cal B}(G)$ are
given, then ${\cal A}_{\mu }(G)= {\cal A}_{\nu }(G)$. From the
definition of measurability of the automorphism the statement of
this lemma follows.

\par {\bf 4. Definitions.} A subgroup $H$ of a group $G$ is called
normal, if its left and right cosets coincide $gH=Hg$ for each $g\in
G$. A group $G$ is called algebraically simple, if it has not a
normal subgroup different from $e$ and $G$, where $e=e_G$ is the
unit element of the group $G$. A topological group $G$ is called
topologically simple, if it has not a normal closed subgroup
different from $e$ and $G$.

\par {\bf 5. Lemma.} {\it  If $f: G\to V$ is a homomorphism of
an algebraically simple group $G$ into a group $V$, then either
$f^{-1}(e_V)=e_G$ or $f^{-1}(e_V)=G$.}
\par {\bf Proof.} The subgroup $J=f^{-1}(e_V)$ is normal in $G$,
since if $f(g)=e_V$, then $f(h^{-1}gh)=
f^{-1}(h)f(g)f(h)=f^{-1}(h)f(h)=f(h^{-1}h)=f(e_G)=e_V$ for each
$h\in G$. Due to the definition of the algebraic simplicity of the
group $G$ either $J=e_G$ or $J=G$.

\par {\bf 6. Corollary.} {\it If an algebraically simple group $G$
has a character $f : G\to S^1$, then either $f^{-1}(1)=e_G$ or
$f^{-1}(1)=G$, where $S^1 := \{ z\in {\bf C}: |z|=1 \} $ is the
multiplicative Abelian group.}

\par {\bf 7. Remark.} If $f^{-1}(e_V)=e_G$ for a homomorphism
$f: G\to V$, then $f$ is bijective. In this case $f(G)$ is
(algebraically) isomorphic with $G$. If moreover $V$ is Abelian, for
example, $V=S^1$, then $G$ is Abelian and can not be simple, besides
the trivial case $G=e$.

\par {\bf 8. Lemma.} {\it If a Hausdorff topological group $G$
is topologically simple and $N$ is a normal subgroup in $G$, then
either $N=e_G$ or ${\bar N} =G$, where ${\bar N}=cl_GN$ is the
closure of $N$ in $G$.}
\par {\bf Proof.} Each topological group is the
uniform space with the entourages of the diagonal of the form $W(U)
:= \{ (g,q)\in G\times G: g^{-1}q\in U \} $, where $U$ is an open
neighborhood of $e$ in $G$ (see Example 8.1.17 in \cite{eng}). If
$g\in {\bar N}$, then there exists a net $\{ g_{\alpha }: \alpha \in
\Lambda \} $ such that $\lim g_{\alpha }=g$, where $\Lambda $ is a
directed set, $g_{\alpha }\in N$ for each $\alpha \in \Lambda $ (see
\S 1.6 and Corollary 8.1.4 in \cite{eng}). Since $h^{-1}g_{\alpha
}h\in N$ for each $h\in G$, then due to the continuity of the
multiplication in $G$ there is satisfied the equality $\lim
h^{-1}g_{\alpha }h=h^{-1}gh$, hence $h^{-1}{\bar N}h={\bar N}$ for
each $h\in G$. Then $\bar N$ is the closed normal subgroup in $G$.
In view of the topological simplicity of the group $G$ either ${\bar
N}=e$ or ${\bar N}=G$. Since $N\subset \bar N$, then in the case of
${\bar N}=e$ we get $N=e$.

\par {\bf 9. Corollary.} {\it If $f: G\to V$ is a continuous
homomorphism of a topologically simple group $G$ in a topological
group $V$, where $G$ and $V$ are Hausdorff, then either
$f^{-1}(e_V)=e_G$ or $f^{-1}(e_V)=G$.}
\par {\bf Proof.} In a Hausdorff topological space
each singleton is closed (see \S 1.5 in \cite{eng}). Since the
homomorphism $f$ is continuous, then $f^{-1} (e_V) =: N$ is closed
in $G$. On the other hand, $N$ is the normal subgroup in $G$. From
Lemma 8 this Corollary follows.

\par {\bf 10. Remark.} Henceforth, Lie groups
$G$ of the smoothness class $C^{\infty }$ over the field of real
numbers or $C^{\omega }$ over local fields, that is finite algebraic
extension of the field of $p$-adic numbers, are considered, where a
smoothness class is for $G$ as the manifold and for the smoothness
of the operation $G\times G\ni (g,q)\mapsto  g^{-1}q\in G$,
$C^{\infty }$ denotes the class of infinitely differentiable
mappings, $C^{\omega }$ denotes the class of locally analytic
mappings. As usually $U(n)$ denotes the unitary group of the unitary
space ${\bf C^n}$, where $n$ is the natural number.

\par {\bf 11. Lemma.} {\it If $G$ is a locally compact Lie group
over $\bf R$ of the dimension $n$, then there exists an open
neighborhood $U$ of the unit element $e$ in $G$, which has a
topological embedding into $(S^1)^n$, as well as an embedding into
$U(n)$ as the local Lie group.}
\par {\bf Proof.} For a locally compact Lie group
$G$ of the smoothness class $C^{\infty }$ over $\bf R$ as it is
well-known there exists a Lie algebra ${\sf g}=T_eG$ and an
exponential mapping $\exp : V_1\to U_1$ of an open neighborhood
$V_1$ of zero in $\sf g$ on an open neighborhood $U_1$ of the unit
element $e$ in $G$, while $\exp $ is the infinite differentiable
diffeomorphism, $U_1$ is a local Lie group (see
\cite{bourgralg,kling,pontr}).
\par As the linear space over $\bf R$ the algebra $\sf g$ has the
dimension $n$, consequently, there exists the embedding of $U_1$
into $\bf R^n$. Choose a compact subset $V$ in $V_1$, such that the
interior $Int (V)$ is the open neighborhood of zero in $\sf g$. Then
$\exp ( V)=: U $ is the compact subset in $U_1$, moreover, $\exp
(Int (V))$ is the local subgroup in $G$. Therefore, $U$ has a
topological embedding into ${\bf R^n}/\bf Z^n$, where the latter
topological space is isomorphic with $(S^1)^n$. Since the Lie
algebra $\sf g$ has the basis of generators $v_1,...,v_m$ of a
dimension not exceeding $n$, then each element in $U$ can be
presented as the finite product of local one-parameter subgroups
$\exp (t_jv_j)$ with $t_j\in (-\epsilon ,\epsilon )$, where
$\epsilon
>0$. For a sufficiently small $\epsilon >0$ each $\{ \exp (t_jv_j):
t_j\in (-\epsilon ,\epsilon ) \} $ has an embedding into $S^1$ as
the Abelian local subgroup. \par In $\sf g$ in the open neighborhood
of zero there is accomplished the Campbell-Hausdorff formula (see
Chapter III in \cite{bourgralg}).
\par Remind that the Campbell-Hausdorff formula for the calculation of
the expression $w=\ln (e^ue^v)$ in the neighborhood of zero $V$ of
the Lie algebra $\sf g$ over ${\bf K}=\bf R$ or a local
non-archimedean field $\bf K$ has the form:
\par $w=\sum_{n=1}^{\infty } n^{-1} \sum_{r+s = n, r\ge 0, s\ge 0}
({{\tilde w}}_{r,s}+{\hat w}_{r,s})$, where \par ${{\tilde w}}_{r,s}
:= \sum_{m\ge 1} (-1)^{m+1} m^{-1} \sum^* ((\prod _{i=1}^{m-1}
(ad\enskip u)^{r_i} (r_i!)^{-1} (ad \enskip v)^{s_i} (s_i!)^{-1})$
\par $ (ad\enskip u)^{r_m} (r_m!)^{-1} )(v),$
\par ${\hat w}_{r,s} :=
\sum_{m\ge 1} (-1)^{m+1} m^{-1} \sum^{**} (\prod _{i=1}^{m-1}
(ad\enskip u)^{r_i} (r_i!)^{-1} (ad \enskip v)^{s_i} (s_i!)^{-1}) (u),$ \\
where $\sum^*$ denotes the sum by $r_1+...+r_m=r,$ $s_1+...+s_{m-1}=
s-1$, $r_1+s_1\ge 1$,...,$r_{m-1}+s_{m-1}\ge 1$, while $\sum^{**}$
denotes the sum by $r_1+...+r_{m-1}= r-1$,
$s_1+...+s_{m-1}=s$, $r_1+s_1\ge 1$,...,$r_{m-1}+s_{m-1}\ge 1$, \\
where the convergence radius of the series depends on $\bf K$ and
the multiplicative norm in it.
\par Each local one-parameter subgroup $\exp (tv)$ with
$t\in (-\epsilon ,\epsilon )$ acts on $\sf g$ and then has the
embedding into $Gl(n,{\bf R})$, therefore, $U$ has the embedding
into $GL(n,{\bf R})$ (see also Theorems 58, 59, 84, 87-90 and
Propositions 42(A,B), 56(A,B,C) in Chapter 10 \S \S 42, 53, 56 and
57 \cite{pontr}).
\par The Lie algebra ${\sf u}(n)$ has the basis
of generators $E_{k,j}-E_{j,k}$, $i(E_{k,j}+E_{j,k})$ with $1\le
k<j\le n$ and $i E_{j,j}$ with $j=1,...,n$, where $i=(-1)^{1/2}$,
$E_{k,j}$ is the $n\times n$ matrix with $1$ on the crossing of the
$k$-th row and $j$-th column and others elements are zero. The Lie
algebra $\sf g$ has the embedding into ${\sf gl}(n,{\bf R})$, while
each generator of the algebra ${\sf gl}(n,{\bf R})$ is the linear
combination of generators of the algebra ${\sf u}(n)$, hence $\sf g$
has the embedding into ${\sf u}(n)$ as the Lie subalgebra. Then
$\exp_u\circ ln_G: U_1\to U(n)$ gives the embedding, where $\exp_u$
is the exponential mapping for the Lie algebra ${\sf u}(n)$ of the
Lie group $U(n)$, $ln_G$ is the logarithmic mapping for $G$ from $U$
into $V$.
\par Certainly, in general an embedding of a local Lie subgroup
may have not an extension over the entire group.

\par {\bf 12. Lemma.} {\it Let $X$ be a Tychonoff (completely
regular) dense in itself topological space with $\sigma $-additive,
$\sigma $-finite non-negative Borel regular measure $\mu $ on a
complete $\sigma $-algebra ${\cal A}_{\mu }$ such that for each
$x\in X$ there exists an open neighborhood $U$ of a finite positive
measure $0<\mu (U)<\infty $, moreover, $\mu $ has not atoms in $X$.
If $f: X\to X$ is the bijective epimorphic mapping such that $card
(f(U)\cap V)\ge {\sf c} := card ({\bf R})$ for each open subsets $U$
and $V$ in $X$, then $f$ and $f^{-1}$ are not $({\cal A}_{\mu
},{\cal B})$-measurable.}
\par {\bf Proof.} Recall, that a measure $\mu $ is called
Borel, if it is defined on the $\sigma $-algebra of all Borel
subsets ${\cal B}(X)$ in $X$, that is the minimal $\sigma $-algebra
generated by the family of all open subsets in $X$. In the given
case the algebra ${\cal A}_{\mu }$ is the minimal $\sigma $-algebra,
produced from ${\cal B}(X)$ and the family of all subsets of $\mu $
measure zero. A measure $\mu $ is called Borel regular, if $\mu (A)
= \sup \{ \mu (C): C\subset A, C \mbox{ is closed } \} $ for each
$A\in {\cal B}(X)$. If $\mu (A)<\infty $, then the transition to the
completion gives $\mu (A) =\inf \{ \mu (V): A\subset V \mbox{ open }
\} $, since the measure $\mu $ is $\sigma $-finite (see Theorem
2.2.2 \cite{federer}).
\par Let $U$ and $V$ be open in $X$, while by the condition of the Lemma
$card (f(U)\cap V)\ge \sf c$, then $f^{-1}(f(U)\cap V)=U\cap
f^{-1}(V)$, since $f$ is the bijective mapping from $X$ onto $X$,
consequently, $card (U\cap f^{-1}(V)\ge \sf c$ for each $U$ and $V$
open in $X$. If $A\in {\cal A}_{\mu }$, $\mu (A)<\infty $, then
there exists a Borel subset $B\in {\cal B}(X)$ such that $A\subset
B$ and $\mu (A)=\mu (B)$. Then for each $\epsilon >0$ there exists
an open subset $W$ in $X$ such that $A\subset W$ and $\mu (A)\le \mu
(W)<\mu (A)+\epsilon $. If $S$ is an arbitrary subset in $X$, then
by $\mu ^*(S)$ there is denoted $\inf \{ \mu (W): S\subset W, W
\mbox{ is open } \} =: \mu ^*(S)$.
\par For arbitrary open subsets $U$ and $V$ of a finite $\mu $
measure take open subsets $U_1$ and $U_2$ in $U$, $V_1$ and $V_2$ in
$V$ such that $0<\mu (U)/2 - \delta < \mu (U_j)<\mu (U)/2+\delta $
and $0< \mu (V)/2 - \delta < \mu (V_j)< \mu (V)/2 + \delta $ for
$j=1, 2$, where $0<\delta <\min (\mu (U), \mu (V))/9$. This is
possible, since $\mu $ is non-negative, has not any atom, while for
each $x\in X$ there exists an open subset $P$ with $x\in P$, $0<\mu
(P)<\infty $. Denote $A:= f^{-1} (U)$, $A_j := f^{-1}(U_j)$. By the
supposition of this lemma $card (A_j\cap V_k)\ge \sf c$ for each $j,
k \in \{ 1, 2 \} $. Suppose that $f$ is a $({\cal A}_{\mu },{\cal
B}(X))$ measurable mapping. Then there would be $A_j\in {\cal
A}_{\mu }$ for $j=1, 2$ and there would exist $B_j\in {\cal B}(X)$
such that $A_j\subset B_j$ and $\mu (B_j)=A_j$, and hence open
subsets would be $W_j$ with $B_j\subset W_j$ and $\mu (B_j)\le \mu
(W_j)< \mu (B_j)+\delta $ for $j=1, 2$. But $A_1\cap A_2=\emptyset
$, consequently, $\mu (A)=\mu (A_1)+\mu (A_2)$. \par If for all open
$U$ there would be $\mu (f^{-1}(U)\cap V)=0$ for each open $V$ in
$X$ with $\mu (V)<\infty $, then in view of the $\sigma $-finiteness
and $\sigma $-additivity of $\mu $ then would be $\mu (X)=0$, that
contradicts the supposition of this lemma, therefore, there can be
chosen open subsets $U$, $U_1$ and $U_2$ such that $\mu (A_j\cap
V)>0$, where $U_1\cup U_2\subset U$. But $card (A_j\cap P_k)\ge \sf
c$ for each $P_k$ - open subset in $V_k$.
\par On the other hand, $\mu ^*(A_j\cap V_k) = \inf \{ \mu (Y): Y
\mbox{ is open }, Y\supset (A_j\cap V_k) \}
>0$, consequently, there exists a countable sequence of open
sets $Y_n\supset (A_j\cap V_k)$ such that $\lim_{n\to \infty } \mu
(Y_n)= \mu ^*(A_j\cap V_k)$. Let $C_n = \bigcap_{s=1}^nY_s$, then
$(A_j\cap V_k)\subset C_{n+l}\subset C_n$ for every $n, l\in \bf N$,
each $C_n$ is open. Moreover, $\lim_{n\to \infty } \mu (C_n)=\mu
^*(A_j\cap V_k)$. Since $card (A_j\cap P)\ge \sf c$ for each $P$
open in $X$, then $ Int (C_n\setminus C_{n+l}) = \emptyset $ for all
$n, l\in \bf N$, where $Int (B)$ denotes the interior of a subset
$B$ in $X$. Then $\mu ^* (A_j\cap V_k) = \mu ^* ([cl_X(A_j\cap
V_k)]\cap V_k)=\mu ([cl_X(A_j\cap V_k)]\cap V_k)$, since
$cl_X(A_j\cap V_k)\in {\cal B}(X)$, while $X$ is the completely
regular space dense in itself, where $cl_X(B)$ denotes the closure
of a subset $B$ in $X$. Thus, $\mu ^*(A_j\cap V_k)\ge \mu
(V)/2-\delta $, since $cl_X(A_j\cap V_k)=cl_X(V_k)\supset V_k$.
Therefore, $\mu (A_1\cap V_1)+\mu (A_1\cap V_2)+\mu (A_2\cap
V_1)+\mu (A_2\cap V_2)\ge 2\mu (V)- 8\delta >10 \mu (V)/9$ and the
contradiction is obtained, since by the construction $A_1\cap
A_2=\emptyset $, $V_1\cap V_2=\emptyset $ and $V_1\cup V_2\subset
V$, consequently, $f$ is not measurable.
\par Applying the above proof to $f^{-1}$ instead of $f$ we get,
that $f^{-1}$ is also nonmeasurable, since $f^{-1}$ satisfies
conditions from the second section of the proof.

\par {\bf 13. Theorem.} {\it Let $G$ be a non-trivial locally
compact Lie group over $\bf R$ or over a non-archimedean local field
$\bf K$. Then the group of its automorphisms $Aut (G)$ has a family
of the cardinality not less than $2^{\sf c}$ of different
non-measurable automorphisms relative to a non-trivial non-negative
Haar measure $\mu $ on $G$, where ${\sf c}:= card ({\bf R})$ denotes
the cardinality of the continuum.}

\par {\bf Proof.} Each locally compact Lie group
$G$ has a finite dimensional Lie algebra $\sf g$ over ${\bf K}=\bf
R$ or over a non-archimedean local field $\bf K$ respectively.
Therefore, both $G$ and $\sf g$ are metrizable, moreover the metric
can be chosen left-invariant (see Theorem 8.3 \cite{hewross}).
\par In the non-archimedean case as a neighborhood $U$
of the unit element $e$ in $G$ there can be chosen a compact clopen
subgroup, since $G$ is totally disconnected (see Theorems 5.13 and
7.7 \cite{hewross}). If ${\bf K}=\bf R$, then we take an open
symmetric $U=U^{-1}$ neighborhood $U$ of the unit element $e$ in
$G$, where $U^{-1} := \{ g^{-1}: g\in U \} $. For a sufficiently
small $U$ there is a bijective exponential mapping $\exp : V\to U$
from the corresponding neighborhood $V$ of zero in $\sf g$ on $U$,
where $\exp $ belongs to the class of smoothness $C^{\infty }$ or
$C^{\omega }$ (see \cite{bourgralg,bourmnog,kling,pontr}). Choose
$U$ sufficiently small, that in it would be satisfied the
Campbell-Hausdorff formula. \par Take a basis of generators
$v_1,...,v_m$ in the Lie algebra $\sf g$ over the field $\bf K$, and
as the linear space over the field $\bf K$ it has the basis $\eta
_1,...,\eta _n$, where $n\ge m$, $v_j=\eta _j$ for $j=1,...,m$, and
$\eta _{m+1},...,\eta _n$ are obtained as finite products
(commutators) $[u,v]$ of basic generators in the Lie algebra $\sf g$
for $n>m$.
\par Consider a set ${\bf K}\setminus {\bf Q}$ of all irrational
elements of the field $\bf K$. The field $\bf K$ is uncountable and
as the linear space over $\bf Q$ it is infinite-dimensional, since
the field of rational numbers $\bf Q$ is countable. For each $b\in
{\bf K}\setminus \bf Q$ there exists an extension ${\bf Q}(b)$ of
the field $\bf Q$ with the help of a number $b$, ${\bf Q}\subset
{\bf Q}(b)$. Remind that a number $a\in {\bf K}\setminus \bf F$ is
called algebraic over the field $\bf F$, if it is a root of a
polynomial with coefficients from $\bf F$. If an element $a$ is not
algebraic over $\bf F$, then it is called transcendental. For a
transcendental element $a$ over a field $\bf F$ the field of
rational quotients with coefficients from the field $\bf F$ is
purely transcendental.
\par A set $\{ a_1,...,a_n \} $ from $\bf K$ is algebraically
independent, if each polynomial $P(x_1,...,x_n)$ with coefficients
from $\bf F$ becoming zero while substitution of $a_1,...,a_n$
instead of $x_1,...,x_n$, is the identically zero polynomial.
Moreover, the field $({\bf F}(a_1,...,a_{n-1}))(a_n)$ is isomorphic
with the field of rational quotients ${\bf F}(a_1,...,a_n)$ of
variables $a_1,...,a_n$ with coefficients from the field $\bf F$. If
a field $\bf F$ is countable, then ${\bf F}(a_j: j\in {\bf N})$ is
also countable, since $\bigcup_{n=1}^{\infty }\aleph _0^n=\aleph
_0$, where $\aleph _0=card ({\bf N})$, $\bf N$ denotes the set  of
natural numbers \cite{eng}.
\par A subset (not necessarily finite) $S$ in ${\bf
K}\setminus \bf F$ is called algebraically independent over $\bf F$,
if each its finite subset is algebraically independent. In the
family of all algebraically independent subsets in $\bf K$ over $\bf
F$ there exists a partial ordering by the inclusion, that makes it
directed. In view of the Kuratowski-Zorn lemma there exists a
maximal relative to such ordering algebraically independent family
$\Psi $ in $\bf K$ over $\bf F$. Each such subset is called the
basis of transcendence of the field $\bf K$ over the field $\bf F$.
In view of Theorem of section 1.1.5 \cite{bacht} the cardinal
numbers of any transcendence bases coincide.
\par Since by the G. Cantor theorem a set $A$ of all
algebraic numbers over a countable field is countable, hence $card
(\Psi )={\sf c}=card ({\bf R})$, in particular, for ${\bf F}=\bf Q$.
Moreover, the Haar measure $\nu $ of algebraic numbers in $\bf K$ is
equal to zero, $\nu (A)=0$, where $\nu $ is non-negative and
non-trivial on $\bf K$. Then ${\bf Q}\subset {\bf Q}(\Psi )\subset
\bf K$, where ${\bf F}(\Psi )$ denotes the purely transcendental
extension of the field $\bf F$, while ${\bf Q}(\Psi )\subset \bf K$
is the algebraic extension (see section 1.1.5 \cite{bacht}).
\par  The Lie algebra $\sf g$ is infinite-dimensional over $\bf Q$ with
an uncountable Hamel basis $\gamma $ over $\bf Q$, that is
$w_1,...,w_k$ are linearly independent over the field of rational
numbers $\bf Q$ for each $w_1,...,w_k\in \gamma $ and $k\in \bf N$,
while each element $w\in \sf g$ is the finite linear combination
$w=c_1w_1+...+c_kw_k$ over the field $\bf Q$ of elements $w_j$ from
$\gamma $ with rational coefficients $c_j\in \bf Q$.

\par In the non-archimedean case $U$ is the subgroup in $G$ (see
above), then put $W=U$. In the case of the Lie group $G$ over the
field of real numbers $\bf R$ the neighborhood $U$ generates the
group $ W:= \bigcup_{n=1}^{\infty }U^n$ containing in it the
connected component $C$ of the unit $e$ (see Theorem 7.4
\cite{hewross}), where $AB:= \{ gh: g\in A, h\in B \} $ for $A,
B\subset G$.

\par Let $g\in G$ be an element of the group $G$,
then elements of the form $g^n$ for $n\in \bf Z$ generate the
commutative subgroup $S(g)$ in $G$. If a subgroup $S(g)$ is finite,
then its order is called the order $ord (g)=k$ of the element $g$,
that is $S(g) = \{ e, g, g^2,...,g^{k-1} \} $ and $g^k=e$. If $S(g)$
is infinite, then it is said, that $g$ is the element of the
infinite order. Let $G_{fin}:= gr ({\cal F})$ be a minimal subgroup
in $G$, generated by all possible finite products of elements from
$\cal F$, where $\cal F$ is the set of all elements of finite orders
in $G$.  Then we consider the subgroup $W_{fin} = G_{fin}\cap W$ in
$G$. Let $U_{fin}=W_{fin}\cap U$. Then $ln (U_{fin})\subset \sf g$,
where $ln : U\to V$ is the logarithmic mapping of a local Lie
subgroup on a neighborhood of zero in the Lie algebra $\sf g$,
corresponding to the Campbell-Hausdorff formula, where coefficients
of the series are rational numbers (see \S 11). Consider the minimal
subalgebra ${\sf g_{fin}}$ over the field of rational numbers $\bf
Q$ generated by  $ln (U_{fin})$ and $v_1,...,v_m$, ${\sf
g_{fin}}\subset \sf g$, since ${\bf Q}\subset \bf K$. \par Consider
one-parameter local subgroups. Each element from the open
neighborhood $U$ of the unit element belongs to a local
one-parameter subgroup $\{ g^t: t\in {\bf K}, |t|<\epsilon \} =
g_{loc}$, where $\epsilon
>0$, which is unique for each non-unit element $g\ne e$
from $U$. Each cyclic commutative group is isomorphic to the group
of roots of $1$ in $\bf C$. Therefore, on each subgroup $g_W :=
\bigcup_{n=1}^{\infty }(g_{loc})^n$ the set of elements of finite
orders has the Haar measure ${\tilde \mu }_g(g_W)=0$ zero, where
${\tilde \mu }_g$ is the Haar measure on $g_W$. Then $\mu
(W_{fin})=0$ and $\nu ^n({\sf g_{fin}})=0$, since the image $\mu
_{ln}$ of the measure $\mu $ on $V$ is equivalent to $\nu ^n|_V$,
where $\nu ^n$ is the Haar measure on $\bf K^n$ as the additive
group, $\mu _{ln}(B) = \mu (\exp (B))$ for each $B\in {\cal B}(V)$.
\par If $g\in G$, $ord (g)=k\in \bf N$ and $h^m=g$, $m\in \bf N$, then
$ord (h)\le mk\in \bf N$. If $g\in U$, $g=e^v$, $v\in V$, $ord
(g)=k\in \bf N$, $h=e^{tv}\in U$ for some $t\in \bf K$, $ord (h)\in
\bf N$, then $t\in \bf Q$.
\par Let ${\bf Q}(A)$ denotes the minimal field, containing all
elements from $A$ and such that ${\bf Q}\subset {\bf Q}(A)\subset
\bf K$. Among subsets $A\subset \Psi $ take such that ${\sf
g}(A)\supset {\sf g_{fin}}$, where ${\sf g}(A)$ is the Lie algebra
over the field ${\bf Q}(A)$ with the Hamel basis $\gamma _A$ as the
${\bf Q}(A)$-linear space, $\gamma _A\subset \gamma $, $\gamma $
denotes the Hamel basis $\sf g$ over $\bf Q$ as the $\bf Q$-linear
subspace (see above). Moreover, it is possible to restrict the
consideration by $A$ such that $(\nu ^n)^*({\sf g}(A))=0$, since
$\nu ^n({\sf g_{fin}})=0$. Among such Lie algebras ${\sf g}(A)$
there exists a minimal due to the Kuratowski-Zorn lemma. Denote it
by ${\sf g}(A_{fin})$, where $A_{fin}\subset \Psi $.
\par If ${\bf F} = {\bf Q}( \{ a_j: j\in \lambda \} ) $, $\lambda
\subset \bf N$, then the field $\bf F$ is countable, since
$\bigcup_{k=1}^{\infty }\aleph _0^k=\aleph _0$. If $B\subset \bf K$
and $\nu (B)=0$, where $\nu $ is a measure equivalent to the Haar
measure on $\bf K$, then $\nu (B+{\bf F})=0$ and $\nu (B{\bf F})=0$
for a countable subfield $\bf F$ in $\bf K$, consequently, $\nu
((B+{\bf F}){\bf F})=0$ and $\nu (\bigcup_{k=1}^{\infty }(B^k+{\bf
F}){\bf F})=0$, since $\nu (B_1B_2)=\int_{\bf K} \chi
_{B_1B_2}(x)\nu (dx)$, where $\chi _B(x)=1$ for $x\in B$, $\chi
_B(x)=0$ for $x\notin B$, $\chi _B$ denotes the characteristic
function of a subset $B$, $B_1+B_2 := \{ x: x=b_1+b_2, b_1\in B_1,
b_2\in B_2 \} $, $B_1B_2 := \{ x: x=b_1b_2, b_1\in B_1, b_2\in B_2
\} $.
\par Since $\nu ^n({\sf g}(A_{fin}))=0$, then $\nu ({\bf Q}(A_{fin
}))=0$, where $\nu _{\bf K}=\nu $ is the Haar measure on $\bf K$ as
the additive group. Consequently, there exists $A_{fin}$ such that
$card (\Phi )=\sf c$, where $\Phi := \Psi \setminus A_{fin}$, since
each element $v$ from $\sf g$ is the finite linear combination over
$\bf Q$ elements from $\gamma $, while in the Campbell-Hausdorff
formula the expansion coefficients are rational (see \S 11).

\par Over each field ${\bf F} = ({\bf Q}(A_{fin}))(b)={\bf Q}(A_{fin}
\cup \{ b \} )$ for $b \in \Phi $ consider a Lie subalgebra ${\sf
g}(A_{fin})_{\bf F}$ generated from the algebra ${\sf g}(A_{fin})$
by extension of the field of scalars ${\bf Q}(A_{fin})$ up to $\bf
F$, ${\sf g}(A_{fin})_{\bf F}\subset {\sf g}$, that is each element
from ${\sf g}(A_{fin})_{\bf F}$ is a finite linear combination over
$\bf F$ elements from ${\sf g}(A_{fin})$.
\par Each element $g\in U$ is a finite product of
elements of local one-parameter subgroups of the form $\exp
(t_jv_j)$, $t_j\in \bf K$. The exponential mapping $\exp : V\to U$
gives an uncountable family of local subgroups $S_{\bf F}=\exp
(V\cap {\sf g}(A_{fin})_{\bf F})$, where ${\bf F}= {\bf Q}(A_{fin
}\cup \{ b \} )$, $b\in \Phi $. All such local subgroups for
different ${\bf F}={\bf Q}(A_{fin}\cup \{ b_j \} )$, $b_1\ne b_2\in
\Phi $, are pairwise isomorphic, since the fields ${\bf
Q}(A_{fin}\cup \{ b_1 \} )$ and ${\bf Q}(A_{fin}\cup \{ b_2 \} )$
are pairwise isomorphic, and in $U$ there is satisfied the
Campbell-Hausdorff formula about the local relation between the
multiplication in the Lie algebra $\sf g$ and the multiplication in
its Lie group $G$ (see Chapters 2 and 3 in \cite{bourgralg}).
Moreover, ${\sf g}(A_{fin})_{{\bf Q}(A_{fin}\cup \{ b_1 \} )}$ is
isomorphic with ${\sf g}(A_{fin})_{{\bf Q}(A_{fin}\cup \{ b_2 \}
)}$. The isomorphism of the fields ${\bf Q}(A_{fin}\cup \{ b_1 \} )$
and ${\bf Q}(A_{fin}\cup \{ b_2 \} )$ is established by the mapping
$\theta =\theta _{b_1,b_2}$, $\theta (b_1)= b_2$ with the identity
mapping on the field ${\bf Q}(A_{fin})$ such that $\theta
(P_k(b_1)/L_s(b_1))= P_k(b_2)/L_s(b_2)$, where $b_1, b_2\in \Phi $,
$P_k$ and $L_s$ are non-zero polynomials of non-negative integer
degrees $k$ and $s$ respectively with coefficients from the field
${\bf Q}(A_{fin})$ (see also \cite{bacht,plotk}).

\par In the non-archimedean case each local subgroup
$S_{\bf F} =: J_{\bf F}$ is the group, since $W=U$ is the group. In
the case of the Lie group $G$ over $\bf R$ take
$W=\bigcup_{n=1}^{\infty }U^n$ (see above). Therefore,
$\bigcup_{n=1}^{\infty }(S_{\bf F})^n =: J_{\bf F}$ is the subgroup
in $G$. Then the group $J_{{\bf Q}(A_{fin}\cup \{ b \} )}$ is
isomorphic with $J_{{\bf Q}(A_{fin}\cup \{ r \} )}$ for each $b\ne
r\in \Phi $.

\par Take a bijective epimorphic mappings $\phi : \Phi \to \Phi $.
It gives an isomorphism $\phi $ of each field ${\bf Q}(A_{fin}\cup
\{ b_1,...,b_z\} )$ onto ${\bf Q}(A_{fin}\cup \{ r_1,...,r_z \} )$,
$\phi (P_k(b_1,...,b_z)/L_s(b_1,...,b_z))=
P_k(r_1,...,r_z)/L_s(r_1,...,r_z)$ for each non-zero polynomials of
integer non-negative degrees $k$ and $s$ of $z\in \bf N$ variables
and with expansion coefficients from ${\bf Q}(A_{fin})$, since $\phi
(x)=x$ for each number $x$ from ${\bf Q}(A_{fin})$, where $r_j=\phi
(b_j)$ for each $j=1,...,z\in \bf N$, $b_1,...,b_z\in \Phi $.
Therefore, it has a natural extension up to an algebraic
automorphism of the field ${\bf Q}(\Psi )=({\bf Q}(A_{fin}))(\Phi )$
onto itself.
\par If $a\in {\bf K}\setminus {\bf Q}(\Psi )$, then the set
$(\Psi ,a)$ is algebraically dependent, that is there exist
$b_1,...,b_k\in \Psi $ such that the family $\{ b_1,...,b_k, a \} $
is algebraically dependent. This means that there exists a
polynomial $F_s(T_1,...,T_k,T_{k+1})$ of a degree $s\in \bf N$ with
rational expansion coefficients such that a substitution of
$b_1,...,b_k$ instead of $T_1,...,T_k$, $a$ instead of $T_{k+1}$ the
polynomial takes the zero value. Then for an automorphism $\phi $ of
the field $\bf K$ a number $\phi (a)$ need to be a root of the
polynomial $F_s(r_1,...,r_k,T_{k+1})$, where $r_j=\phi (b_j)$ for
every $j=1,...,k$. If $a$ is an algebraic number over $\bf Q$, then
$k=0$ and there can be taken $\phi (a)=a$. \par For each algebraic
number $a$ from ${\bf K}$ over the field ${\bf Q}$ put $\phi (a)=a$.
Since the field $\bf Q$ is everywhere dense in the field $\bf R$ or
$\bf Q_p$ respectively, while the fields $\bf R$ and $\bf Q_p$ are
complete as the normed spaces relative to their multiplicative
norms, then $\bf Q$ is everywhere dense in ${\bf Q}(\Psi )$ relative
to the norm inherited from the field $\bf K$. Therefore, every
algebraic number from the field $\bf K$ over the field ${\bf Q}(\Psi
)$ is the limit of a converging sequence of algebraic numbers from
$\bf K$ over the field $\bf Q$. The field $\bf C$ is algebraically
closed and is obtained by the way of the extension of the field $\bf
R$ with the help of the root of the polynomial $x^2+1=0$, which is
invariant relative to the automorphism $\phi $.
\par In the case of the non-archimedean local field $\bf K$ the residue
class field $B({\bf K},0,1)/B({\bf K},0,|\pi |)$ is the finite field
${\bf F}_q$ with the number $q=p^y$ of its elements, where $y\in \bf
N$, $\pi \in \bf K$, $|\pi |= \max \{ |x|: x\in {\bf K}, |x|<1 \} $,
$B({\bf K},x_0,r) := \{ x\in {\bf K}: |x-x_0|\le r \} $ is the ball
of radius $r>0$ in $\bf K$, containing $x_0$. At the same time the
residue class field $B({\bf Q_p},0,1)/B({\bf Q_p},0,1/p)$ is
composed of $p$ elements, where $|p|=1/p$ for $p\in \bf Q_p$. The
field of $p$-adic numbers has the normalization group $\Gamma _{\bf
Q_p} := \{ |x|: x\ne 0, x\in {\bf Q_p} \} = \{ p^k: k\in {\bf Z} \}
$. While $\Gamma _{\bf K} = \{ p^{k/l}: k\in {\bf Z} \} $ for some
$0<l\in \bf N$. That is $\bf K$ is obtained by a finite algebraic
extension of the field of $p$-adic numbers by adding roots of
polynomials with expansion coefficients from the field of rational
numbers $\bf Q$ (see also Theorem 7 and Proposition 5 in \S I.4
\cite{weil}).
\par Therefore, for each algebraic number $a\in {\bf
K}\setminus {\bf Q}(\Psi )$ there exists a polynomial
$F_s(b_1,...,b_v,X)$ with $b_1,...,b_v\in \Psi $ of degree $s\in \bf
N$ with rational expansion coefficients such that
$F_s(b_1,...,b_v,a)=0$, consequently, it can be taken $\phi (a)=c$,
where $F_s(r_1,...,r_v,c)=0$, $r_j=\phi (b_j)$, $c\in \bf K$, since
$\phi (q)=q$ for each rational number $q\in \bf Q$, while ${\bf
Q}(b_1,...,b_v)$ and ${\bf Q}(r_1,...,r_v)\subset {\bf Q}_p$. Thus,
the automorphism $\phi $ has an extension up to an automorphism
$\phi : {\bf K}\to \bf K$ as in the case ${\bf K}=\bf R$, as well as
for the local field $\bf K$. \par Since $card (\Phi )=\sf c$, $card
{\bf Q}(\Phi )=\sf c$ and ${\bf Q}(\Phi )$ is everywhere dense in
$\bf K$, then there exist bijective surjective mappings $\phi : \Phi
\to \Phi $ generating automorphisms of the field $\phi : {\bf K}\to
\bf K$ as above such that for each open subsets $S$ and $T$ in $\bf
K$ there is satisfied the relation for the cardinality $card (\phi
(S)\cap T)=\sf c$. Indeed, the set $\Phi $ can be described in the
form of the disjoint union of subsets $\Lambda _a$, $a\in E$, $card
(E)=\sf c$, $card (\Lambda _a)=\sf c$ for each $a$, $\bigcup_{a\in
E}\Lambda _a= \Phi $, $\Lambda _a\cap \Lambda _b=\emptyset $ for
each $a\ne b$. At the same time $card ({\bf Q}(\Lambda _a)\cap
T)=\sf c$ for each $T$ open in $\bf K$. Since $card (\Lambda
_a)=card (\Lambda _b)$, then there exists a bijection $\phi ^a_b:
\Lambda _a\to \Lambda _b$ from $\Lambda _a$ onto $\Lambda _b$ for
each $a, b$.
\par Take a bijective mapping $\eta : E\to E$ from $E$ onto
$E$ such that $\eta (a)\ne b$ for each $a$. Then the combination of
mappings $\{ \phi ^a_{\eta (a)}: a\in E \} $ generates the bijective
mappings $\phi : \Phi \to \Phi $, putting $\phi (a)=a$ on $A_{fin}$,
we get the bijective mapping $\phi : \Psi \to \Psi $ from $\Psi $ on
$\Psi $. In view of the proof given above it has an extension up to
an automorphism of the field $\phi : {\bf K}\to \bf K$. Therefore,
$card (\phi (S)\cap T)=\sf c$ for each $S$ and $T$ open in $\bf K$,
since $card ({\bf Q}(\Lambda _a)\cap T) = \sf c$ for each $T$ open
in $\bf K$. It is not difficult to mention, that the family of such
different algebraic automorphisms of the field $\bf K$ has the
cardinality $2^{\sf c}$, since $card ({\sf c}^{\sf c})=2^{\sf c}$
\cite{eng}. In view of Lemma 12 they are $({\cal A}_{\nu },{\cal
B}({\bf K}))$-non-measurable, where $\nu $ - is the non-negative
non-trivial Haar measure on $\bf K$ (see also \cite{bourhm,weil}).
\par On the other hand, every automorphism $\phi $ of the field $\bf K$
generates an automorphism $f$ of the group $G_0 := gr(J_{{\bf Q}
(A_{fin})}\cup [\bigcup_{b\in \Phi }J_{{\bf Q}(A_{fin}\cup \{ b \}
)}])$, where $gr (B)$ is a minimal algebraic subgroup in $G$
generated by finite products $g_1^{a_1}...g_r^{a_r}$ of all elements
$g_j\in B$, $a_j\in \bf Z$, $j=1,...,r\in \bf N$. The automorphism
$f$ is produced with the help of the isomorphisms $f: J_{{\bf
Q}(A_{fin}\cup \{ b \} )} \to J_{{\bf Q}(A_{fin}\cup \{ \phi (b) \}
)}$ for each $b\in \Phi $, while the restriction $f| J_{{\bf
Q}(A_{fin})}$ take, for example, the identity mapping. The mapping
$f$ has an extension from $A := J_{{\bf Q}(A_{fin})}\cup
[\bigcup_{b\in \Phi }J_{{\bf Q}(A_{fin} \cup \{ b \} )}]$ on $G_0$
by taking of all possible finite products of initial elements from
$A$. For each such subgroup $G_0$ of the group $G$ the automorphism
has an extension on $G$, since every automorphism $q$ of a subgroup
$Y$ of the group $G$ can be extended up to an automorphism of some
group $H$ containing in itself $G$, but $q(G)$ is contained in $H$
and is isomorphic with $G$ (see \cite{focus,neumann}).
\par By the construction of the group $G_0$ it is everywhere dense
in the group $W$ in the topology inherited from $G$. Then due to
Lemma 12 the automorphism $f$ is $({\cal A}_{\mu },{\cal B}(G))$
non-measurable. This also follows from the fact that the exponential
mapping $\exp $ from the neighborhood of zero $V_0$ in the Lie
algebra $\sf g$ on the neighborhood $U_e$ of the unit element in $G$
induces the image of the measure $\nu ^n_{\exp }$ on $U_e$, where
$n$ is the dimension of $\sf g$ as the linear space over the field
$\bf K$, $\nu ^n$ is the Haar measure on $\bf K^n$ as the additive
group. At the same time the measure $\nu ^n_{\exp }$ is equvalent to
the restriction of the Haar measure $\mu $ on $U_e$ (see Definitions
1 above).
\par From the proof of Theorem 13 it follows the following.

\par {\bf 14. Corollary.} {\it The family $\Upsilon $ of non-measurable
automorphisms from Theorem 13 has a subfamily  $\Omega $ of the
cardinality $card (\Omega ) \ge 2^{\sf c}$ such that every $f\in
\Omega $ after a restriction on each one-parameter subgroup over the
field $\bf K$ in $G$ is non-measurable relative to the corresponding
Haar measure on the subgroup.}

\par {\bf 15. Theorem.} {\it Let $\sf g$ be a non-trivial
Lie algebra finite-dimensional over the field $\bf K$ with a measure
$\mu $ equal to the non-trivial non-negative Haar measure on the
additive group for $\sf g$. Then the algebra $\sf g$ has $2^{\bf c}$
non-measurable automorphisms.}
\par {\bf Proof.} Take any algebraic automorphism $\phi $ of the
field $\bf K$ from the proof of Theorem 13. Since $\sf g$ is
finite-dimensional over the field $\bf K$, then the non-negative
non-trivial Haar measure $\mu $ on $\sf g$ as the additive group for
$\sf g$ is equivalent to the measure $\nu ^n$ (see Definitions 1).
The automorphism $\phi $ has the extension up to an automorphism of
the algebra: $\phi (a_jv_j)=\phi (a_j)v_j$, $\phi
(a_1v_1+....+a_mv_m)= \phi (a_1)v_1+...+\phi (a_m)v_m$, $\phi
([a_kv_k,a_jv_j]) = [\phi (a_k)v_k, \phi (a_j)v_j]$ for each $a_j\in
\bf K$, $k, j=1,...,m$, where $v_1,...,v_m$ is the basis of
generators in $\sf g$  (see Definitions 2.1). Therefore, in view of
Lemma 12 the automorphism $\phi $ is $({\cal A}_{\mu },{\cal B}({\sf
g}))$-non-measurable, relative to the additive group of the algebra
$\sf  g$. The family of such different automorphisms of the algebra
$\sf g$ has the cardinality $2^{\sf c}$.

\par {\bf 16. Theorem.} {\it Let $G$ be a locally compact Hausdorff
group with a countable base of neighborhoods of the unit element $e$
and $f: G\to G$ be its automorphism non-measurable relative to the
non-trivial non-negative left- (or right-)invariant Haar measure
$\mu $ on $G$. Then $G$ has a topologically irreducible unitary
representation, which is not weakly measurable.}
\par {\bf Proof.} Since $G$ has a non-measurable
automorphism, then it is non-discrete, since all algebraic
automorphisms are continuous relative to the discrete topology,
since in discrete topology every point is an open subset. On the
other hand, every $T_0$ topological group is completely regular (see
Theorem 8.4 \cite{hewross}). In view of Theorem 5.8 \cite{hewross} a
subgroup of a topological group is discrete if and only if it
contains an isolated point. Therefore, the group $G$ is dense in
itself, that is every its point $q$ is a limit of a convergent net,
contained in the punctured open neighborhood $U\setminus \{ q \} $
of a point $q$. In view of Lemma 5.28 \cite{hewross} locally
countably compact regular topological space $Y$ can not be presented
as a countable union of closed subsets with the empty interior.
Thus, every open subset $U$ in $G$ is uncountable, $card (U)\ge \sf
c$ (see Remark (4.26) \cite{hewross}). In view of this the Haar
measure $\mu $ on $G$ has not atoms. By its construction the Haar
measure is Borel regular.
\par Recall, that the unitary representation is a homomorphism
$T: G\to U(X)$, where $U(X)$ is the unitary group of the Hilbert or
the unitary space $X$ over the field of complex numbers $\bf C$. It
is called topologically irreducible, if in $X$ there does not exist
any closed invariant subspace relative to the family of unitary
operators $\{ T_g: g\in G \} $, besides $\{ 0 \} $ and the entire
$X$.
\par A unitary representation is called weakly measurable, if
the functions $(y,T_gz)$ are $({\cal A}_{\mu }, {\cal B}({\bf
C}))$-measurable for each given vectors $y, z\in X$, where $g\in G$.
On every locally compact group there exists a non-trivial
non-negative left-(or right-)invariant Haar measure (see \S 27
\cite{nai} and \cite{bourhm,hewross}).
\par Then the space $L^1(G,{\cal A},\mu ,{\bf C})$ is supplied with
the ring structure with the usual addition of functions and
convolutions of functions as multiplication in the ring, which is
called the group ring  (see \S 28 \cite{nai}). A ring $R$ is called
normed, if it is the normed space, for each $x, z\in R$: $|xz|\le
|x| |z|$, and if in $R$ there is the unit $e$, then $|e|=1$. The
complete normed ring is called the Banach ring. A ring $R$ is called
symmetric, if it is supplied with the involution $x\mapsto x^*$,
which maps from $R$ into $R$. In particular, the ring $L^1(G,{\cal
A},\mu ,{\bf C})$ is symmetric. \par Therefore, a representation of
the group generates a representation of the group ring in the ring
(algebra in more modern terminology) $L(X)$ of bounded linear
operators from $X$ into $X$. The adjoining of the unit to the ring
$L^1(G,{\cal A},\mu ,{\bf C})$ gives the group ring with the unit,
denote it by $R(G)$. In view of Theorem  29.1 \cite{nai} to each
representation $x\mapsto A_x$ of the group ring $R(G)$ not
containing a degenerate representation, there corresponds a
continuous unitary representation $g\mapsto T_g$ of the group $G$.
Vice versa, to every weakly continuous unitary representation
$g\mapsto T_g$ of the group $G$ there corresponds a representation
$x\mapsto A_x$ of its group ring $R(G)$, which does not contain a
degenerate representation. These representations are related with
each other by the formula: $A_{be+f} = b1 + \int_G f(g)T_g \mu (dg)$
for each $b\in \bf C$, $f\in L^1(G,{\cal A},\mu ,{\bf C})$.
\par A linear functional $f$ is called positive, if
$f(x^*x)\ge 0$ for each $x\in R$. If $f_1$ and $f$ are positive
functionals, then $f_1$ is called subordinated to the functional
$f$, which is denoted by $f_1<f$, if there exists a number $b$ such
that $bf-f_1$ is a positive functional in $R$. A functional $f_1$ in
a symmetric ring is called subordinated to a given positive
functional $f$, if $f_1$ is a linear combination with complex
coefficients of positive functionals subordinated to the functional
$f$. A positive functional $f$ is called indecomposable, if each
functional $f_1$ subordinated to the functional $f$ is its multiple,
that is $f_1=bf$, where $b\in \bf C$.  A representation $x\mapsto
A_x$ is called cyclic, if there exists a vector $y_0\in X$ such that
$ \{ A_xy_0: x\in R \} $ is everywhere dense in $X$.
\par In view of Theorem 19.3.1 \cite{nai} a cyclic representation
of a Banach symmetric ring $R$: $x\mapsto A_x$ is irreducible if and
only if, each defining it positive functional $f(x)=(A_xy_0,y_0)$ is
indecomposable. Let $S$ be a set of all positive functionals on $R$
such that $f(e)=1$. In view of Proposition 19.4.1 \cite{nai} a
positive functional $f$, satisfying condition $f(e)=1$ is
indecomposable if and only if it is an extremal point in the set
$S$.
\par If $H$ is a subgroup in $G$ and $M$ is a topologically
irreducible unitary representation of $H$ in the Hilbert space $Y$,
$Y\ni y_0\ne 0$, $\| y_0 \| =1$, $t(h):= (M_hy_0,y_0)$, $h\in H$,
then the function $t$ is positive definite. Consider a convex subset
$W$ of all positive definite functions on $G$, coinciding with $t$
on $H$. This set is non-void, since it contains the function equal
to $t$ on $H$ and $0$ on $G\setminus H$. If $v\in W$, then
$v(e)=(T_ey_0,y_0)=1$, consequently, $|v(g)|\le 1$ for every $g\in
G$. Therefore, $W$ is compact in the topology of pointwise
convergence. In view of the Krein-Milman theorem $W$ contains an
extremal point $r$ (see Theorem 3.9.1 \cite{nai}), its restriction
is $r|_H=f$. Consequently, $r$ is indecomposable and to it a
topologically irreducible unitary representation of the group $G$
corresponds (see \S 2.1 \cite{bicht}, \cite{nai,fell,hewross}).

\par Consider the Hilbert space $X := L^2(G,{\cal A},\mu ,{\bf
C})$ of the equivalence classes of all functions $v: G\to \bf C$
with square integrable its module on $G$ relative to the measure
$\mu $. Then  there exists a strongly continuous unitary regular
representation $T: G\to U(X)$. The strong continuity means, that
$T_gz$ is the continuous mapping from $G$ into $X$ for each $z\in
X$, where $X$ is supplied wit the standard norm associated with the
scalar product: $\| z \| ^2 = (z,z)$. In the case of the
left-invariant Haar measure it is given by the formula $T_gv(h) :=
v(g^{-1}h)$ for each $g, h\in G$, $v\in X$, where $\mu (gJ)=\mu (J)$
for every $g\in G$ and each $\mu $-measurable subset of finite
measure, $J\in {\cal A}$.
\par Evidently, if a mapping (in the given case a functional
$(T_gx,y)$ of the representation $T$) is non-measurable on an open
subset $W$ in $G$, then it is non-measurable on $G$. For the proof
of non-measurability it is sufficient to take a clopen subgroup $W$
in $G$, which is compact, if $G$ is totally disconnected, or
$W=\bigcup_{n=1}^{\infty }U^n$ for a locally connected $G$ which is
not totally disconnected, where $U$ is a symmetric open neighborhood
of $e$ in $G$, $0<\mu (U)<\infty $ (see Theorems 7.7 and 5.7
\cite{hewross}). Since $e$ has the countable base of neighborhoods
of $e$, then the spaces $L^p(W,{\cal A},\mu ,{\bf C})$ for $1\le
p<\infty $ are separable. Then the proof reduces to the
consideration of the subgroup $W$. Denote $W$ by $G$.

\par The regular representation $T$ is injective, $T_g\ne T_h$
for each $g\ne h\in G$. In view of Theorem 41.4.3 \cite{nai} it can
be decomposed into the direct integral of topologically irreducible
unitary representations $T= \bigoplus \int_S T^s\lambda (ds)$, where
$S$ is a compact (bi-compact in old terminology) Hausdorff
topological space, where $\lambda $ is a $\sigma $-additive measure
on ${\cal B}(S)$ (see also \cite{eng,nai,hewross}). At the same time
the representation $T^s$ is strongly continuous. Therefore,
$q^{-1}({\cal B}(G))\subset {\cal B}(X)$ for every function $q(g) :=
T^s_gz$ for marked $z\in X$, $s\in S$, where $q: G\to X$.
\par For each non-zero vector $z\in X^s$, $z\ne 0$, the closure
of the linear span over the field of complex numbers $\bf C$ of all
vectors $T^s_gz$ coincides with $X^s$, where $X^s$ is an invariant
closed subspace in $X$ relative to the unitary representation $T^s$.
Then the set of functions $\{ q_{x,z,s}: G\to {\bf C}; x\in X^s \} $
separates points in $X^s$, where $q_{x,z,s}(g) := (x,T^s_gz)$, $x, z
\in X^s$. At the same time $q_{x,z,s}\circ f(g)= (x,T^s_{f(g)}z)$.
If $A\subset \bf C$, then $(q_{x,z,s}\circ f)^{-1}(A)= f^{-1}(
q_{x,z,s}^{-1}(A))$. If $A$ is open in $\bf C$, then
$q_{x,z,s}^{-1}(A)$ is open in $G$. Since $f$ is non-measurable on
$G$, also $G\ni h\mapsto gh\in G$ is the continuous mapping from $G$
onto $G$, $gU$ is open for each $g\in G$ and open $U$ in $G$, then
the restriction $f|_U$ is non-measurable for each $U$ open in $G$.
\par Consider an algebraic homomorphism $T\circ f: G \to U(X)$.
If every $T^s\circ f$ would be weakly $({\cal A}_{\mu },{\cal
B}({\bf C}))$-measurable, then $T\circ f$ also would be weakly
$({\cal A}_{\mu },{\cal B}({\bf C}))$-measurable. But in view of
strong continuity of $T$ and non-measurability of the automorphism
$f$ the composition $T\circ f$ is not weakly $({\cal A}_{\mu },{\cal
B}({\bf C}))$-measurable (see also Lemma 12 and Theorem 13 above).
\par {\bf 17. Corollary.} {\it If $G$ is a non-trivial locally
compact Lie group, then it has not less than $2^{\sf c}$ weakly
non-measurable relative to a non-trivial non-negative Haar measure
$\mu $ topologically irreducible unitary representations.}
\par {\bf 18. Remark.} From the proofs of Theorems 13 and 16 it follows,
that the property of $\mu $-non-measurability of an automorphism $f$
or of a unitary representation is local: if $f$ is $\mu
$-non-measurable for the restriction on an open subset $W$ in $G$,
then it is $\mu $-non-measurable on $G$. \par For a totally
disconnected locally compact group it can be taken as $W$ a clopen
compact subgroup in $G$ (see Theorem 7.7 \cite{hewross}). In the
case of locally connected which is not totally disconnected group
$G$ it can be taken a clopen subgroup $W=\bigcup_{n=1}^{\infty
}U^n$, where $0<\mu (U)<\infty $, $U$ is a connected symmetric open
neighborhood of $e$ in $G$ (see Theorem 5.7 \cite{hewross}). Then on
$W$ there exists a probability measure equivalent with $\mu $ (see
also Lemmas 2 and 3). In the case of a locally compact Lie group
over $\bf R$ it can be taken on $U$ also a measure equivalent with
the measure $\lambda _{\theta ^{-1}}$, where $\theta : U\to (S^1)^n$
is the topological embedding from Lemma 11, while $\lambda $ is the
Haar measure on $(S^1)^n$, $\lambda _{\theta ^{-1}}(Y) := \lambda
(\theta (Y))$ for each $Y\in {\cal B}(U)$.
\par In Theorem 16 and Lemma 17 there are irreducible unitary
representations, since in general one can not restrict on characters
because of Lemmas 5, 8 and Corollaries 6,9.

\par {\bf 19.} Let $G$ be a $C^{\infty }$ or $C^{\omega }$
non-trivial Lie group over the field ${\bf K}=\bf R$ or a
non-archimedean local field $\bf K$, moreover, $G$ is complete as an
uniform space. Suppose that $\mu $ is a $\sigma $-additive $\sigma
$-finite non-negative non-trivial Borel regular measure on ${\cal
B}(G)$ such that for each $g\in G$ there exists an open neighborhood
$U$, $g\in U$, with $0<\mu (U)<\infty $, moreover, $\mu $ has not
any atom. Let $\exp : V\to U$ be the exponential mapping for $G$ as
the smooth manifold over $\bf K$ from the open neighborhood $V$ of
zero in $T_eG$ onto an open neighborhood $U$ of the unit element $e$
in $G$. Let also a measure $\nu $ on $T_eG$ be such that $\nu
(J)=\mu (\exp (J))$ for each Borel subset $J\in {\cal B}(T_eG)$,
where $T_eG$ is the linear space of separable type over the field
$\bf K$.

\par Suppose that $G$ has an open neighborhood $U$ of the
unit element $e$ with $0<\mu (U)<\infty $, $\exp : V\to U$ for which
a set of elements $S := \{ g\in U: g \mbox{ belongs to a local}$
$\mbox{ one-parameter}$ \\ $\mbox{ subgroup in }$ $U, \mbox{ over
the field } {\bf K},$ $\mbox{ moreover it is unique}$ $\mbox{for a
given } g \} $ has a positive outer measure, $\mu ^*(S)>0$, where
local one-parameter subgroups have the form $\{ \exp (xv): x\in {\bf
K}, |x|<\epsilon \} \subset U$, $v\in T_eG$, $\epsilon
>0$. Suppose that the restriction of $\ln $ on $U$ corresponds to the
Campbell-Hausdorff formula, where $\ln $ is the inverse mapping to
$\exp $.
\par Let $\pi _v: T_eG\to {\bf K}v$ be a linear over $\bf K$ projection
operator, moreover, $\nu _v$ is equivalent to the  Haar measure on
$\bf K$, where $v$ is a non-zero vector of the tangent space $v\in
T_eG$, ${\bf K}v$ is the one-dimensional over $\bf K$ subspace in
$T_eG$ containing a vector $v$, $\nu _v(J)=\nu (\pi _v^{-1}(J))$ for
each $J\in {\cal B}({\bf K})$.

\par {\bf Theorem.} {\it Then such group $G$ has a family
$\Upsilon $ of $({\cal A}_{\mu },{\cal B}(G))$-non-measurable
automorphisms of the cardinality not less than $2^{\sf c}$, $card
(\Upsilon ) \ge 2^{\sf c}$. Moreover, $\Upsilon $ has a subfamily
$\cal P$ of automorphisms $f$, restrictions of which on
one-parameter over $\bf K$ local subgroups $\{ \exp (xv):
|x|<\epsilon \} $ in $S$ are non-measurable relative to the
corresponding Haar measure on $\{ \exp (xv): |x|<\epsilon \} $.}
\par {\bf Proof} is proceeded by the generalization of the proof
of Theorem 13. For this consider a minimal subgroup $G_S$ in $G$,
generated by elements in $S$. The image $\nu $ on $V$ of the measure
$\mu $ with the help of the logarithmic mapping $\ln $, $\nu (B)=
\mu (\exp (B))$ for every $B\in {\cal B}(V)$, has the extension up
to a $\sigma $-additive finite measure on $\sf g$, $\nu (B) :=
\sum_{j=1}^{\infty } \nu ((B-h_j)\cap V)/2^j$ for every $B\in {\cal
B}({\sf g})$, where $\{ (V+h_j): j\in {\bf N}, h_j\in {\sf g} \} $
is the covering of $\sf g$.  Therefore, $\nu $ has $\sigma
$-additive projections $\nu _{\omega ({\bf K^n})}$ on $\omega ({\bf
K^n})$ for each embedding $\omega : {\bf K^n} \hookrightarrow \sf g$
as the $\bf K$-linear space, $\nu _{\omega ({\bf K^n})}(B) = \nu
(\pi ^{-1}(B))$ for each $B\in {\cal B}(\omega ({\bf K^n}))$, where
$\pi : {\sf g}\to \omega ({\bf K^n})$ is the projection.
\par The operator $\pi $ is $\bf K$-linear and it exists, since
$\omega ({\bf K^n})$ is finite-dimensional over $\bf K$, the field
$\bf K$ is locally compact (see Theorems 5.13 and 5.16 \cite{rooij}
and \cite{nari}). The image $\nu _{\omega ({\bf K})}|_{V\cap \omega
({\bf K})}$ with the help of $\exp $ generates the measure on a
local one-parameter subgroup in $S$, $\mu _g(B)= \nu _{\omega ({\bf
K})}(\ln (B))$ for every $B\in {\cal B}(g_W\cap U)$, which has the
extension up to the $\sigma $-additive measure $\mu _g$ on $g_W$.
Therefore, the set of elements of finite orders from $g_W$ has $\mu
_g$-measure zero. Then as in \S 13 $\mu ^*((G_S)_{fin}\cap S)=0$,
since the measure $\mu $ is Borel regular, $\sigma $-additive and it
has not atoms. Since the  measure $\mu $ is $\sigma $-finite, then
also $\mu ^*((G_S)_{fin})=0$.
\par Consider the Lie algebra ${\sf g}$ generated by $\ln (S)$ over
the field $\bf K$. As the linear space $\sf g$ has the separable
type over $\bf K$, that is there exists a countable family $\rho $
of $\bf K$-linearly independent vectors the linear span of which
$span_{\bf K}\rho $ is everywhere dense in $\sf g$. Then the minimal
Lie algebra over $\bf K$ $Lie alg (span_{\bf K}\rho )$ containing
$span_{\bf K}\rho $ is everywhere dense in $\sf g$.
\par Construct $\sf g_{fin}$ as the minimal algebra over the field
of rational numbers $\bf Q$ generated by $\ln ((G_S)_{fin}\cap S)=0$
and $\rho $. For every minimal subalgebra ${\sf g}(\{ v_j: j\in
\lambda \} )$ over $\bf K$ with generators $v_j\in \rho $, $\lambda
\subset \bf N$, there exists $A\subset \Psi $ such that ${\sf g}(\{
v_j: j\in {\bf N} \} )\cap {\sf g_{fin}}\subset {\sf g}(A)$, where
${\sf g}(A)$ is the minimal Lie algebra over the field ${\bf Q}(A)$,
satisfying this inclusion, $\gamma _A$ denotes the Hamel basis ${\sf
g}(A)$ over ${\bf Q}(A)$, since every $v$ from $\sf g$ is the finite
linear combination over $\bf Q$ of elements from $\gamma $, also
$card (\bigcup_{k=1}^{\infty }(\gamma \aleph _0)^k)=card (\gamma
\aleph _0)=card (\gamma )\ge \sf c$, while expansion coefficients in
the Campbell-Hausdorff formula are rational numbers, $\gamma $
denotes the Hamel basis of $\sf g$ over $\bf Q$ as the $\bf
Q$-linear space, $card (\Psi \setminus A)=\sf c$.
\par If $B\subset \bf K$ and $\nu (B)=0$, where $\nu $ is a measure
equivalent with the Haar measure on $\bf K$, then $\nu (B+{\bf
F})=0$ and $\nu (B{\bf F})=0$ for a countable subfield $\bf F$ in
$\bf K$, consequently, $\nu ((B+{\bf F}){\bf F})=0$ and $\nu
(\bigcup_{k=1}^{\infty }(B^k+{\bf F}){\bf F})=0$, since $\nu
(B_1B_2)=\int_{\bf K} \chi _{B_1B_2}(x)\nu (dx)$, where $\chi
_B(x)=1$ for $x\in B$, $\chi _B(x)=0$ for $x\notin B$, $\chi _B$ is
the characteristic function of a subset $B$, $B_1+B_2 := \{ x:
x=b_1+b_2, b_1\in B_1, b_2\in B_2 \} $, $B_1B_2 := \{ x: x=b_1b_2,
b_1\in B_1, b_2\in B_2 \} $.
\par Since $\mu ^*((G_S)_{fin})=0$ and $\nu _{\sf g}^*({\sf g_{fin}})=0$,
$card (\rho )\le \aleph _0$, then there exists $A_{fin}\subset \Psi
$ such that ${\sf g}(A_{fin})\supset {\sf g}_{fin}$ and $card (\Psi
\setminus A_{fin})=\sf c$, since $\beta \aleph _0=\beta $ for every
cardinal number $\beta \ge \aleph _0$, moreover, $\nu _{\omega ({\bf
K})}({\bf Q}(A_{fin}))=0$, since $\nu _{\omega ({\bf K})}(B)>0$ for
each $B$ open in $\omega ({\bf K})$, $\nu _{\omega ({\bf K})}$ is
equivalent to the Haar measure on $\bf K$ due to the conditions from
\S 19.
\par Take $\Phi = \Psi \setminus A_{fin}$ and an algebraic
automorphism $\phi $ of the field $\bf K$ from \S 13.  Analogously
to \S 13 we construct an automorphism $f$ of the group $G_S$, also
it has an extension up to an automorphism of the group $G$ (see
\cite{neumann,focus}). In view of Lemma 12 and $\mu ^*(S)>0$ it
follows, that $f$ is $({\cal A}_{\mu },{\cal B}(G))$-non-measurable.
Non-measurability of restrictions of $f$ on one-parameter subgroups
relative to the Haar measure on them follows from the properties of
$\phi $ as in \S 13.

\par {\bf 20. Theorem.} {\it Let $G$ be an infinite topological
dense in itself Hausdorff group with a non-negative non-trivial
Borel regular measure $\mu $ on ${\cal B}(G)$ having not any atom,
also for each $g\in G$ there exists an open neighborhood $U$ such
that $0<\mu (U)<\infty $, moreover, $G$ is complete as the uniform
space, $card (U)\ge \sf c$ for each open $U$ in $G$. Then there
exists a family $\Upsilon $ of $({\cal A}_{\mu },{\cal
B}(G))$-non-measurable different automorphisms of the group $G$ of
the cardinality $card (\Upsilon )\ge 2^{\sf c}$.}
\par {\bf Proof.} Take an open symmetric neighborhood
$U$ of the unit element $e$ in $G$ such that $0<\mu (U)<\infty $,
then $\mu $ is the $\sigma $-finite measure on
$W=\bigcup_{n=1}^{\infty }U^n$ (see Theorem 7.4 \cite{hewross}). In
view of the fact that  $G$ is infinite and dense in itself, then it
is non-discrete (see Theorem 5.8 \cite{hewross}). If $\{ g_n: n\in
\alpha \} $ is a net converging to $e$, $\alpha $ is an ordinal,
$card (\alpha )\ge \aleph _0$, then $\{ gg_n: n\in \alpha \} $ is
the net converging to $g$ for each $g\in G$.
\par The topological Hausdorff group $G$ can be supplied with the
left-invariant uniformity giving the initial topology on $G$, while
the left-invariant uniformity can be produced with the help of the
family $\{ \eta _x: x\in M \} $ left-invariant pseudo-metrics $\eta
_x(g,h)=\eta _x(h^{-1}g,e)$ for each $g, h\in G$, then $\eta _x
(g^a,g^b)=\eta _x(g^{a-b},e)$ for each $g\in G$ and each $a, b\in
\bf Z$, where $M$ is some set (see Chapter 8 in \cite{eng}).
\par If $g\in G$, $ord (g)=k<\infty $, then $ord (g^r)\le k$
for each $r\in \bf Z$, in particular, for mutually prime numbers $r$
and $k$ there is satisfied $ord (g)=ord (g^r)$. If $ord (g)=\omega
_0$, then $ord (g^r)=\omega _0$ for each integer non-zero $r$, where
$\omega _0$ denotes the initial ordinal of the cardinality $\aleph
_0$.
\par Since the measure $\mu $ is non-trivial, non-negative and for
$e\in G$ there exists an open neighborhood such that $0<\mu
(U)<\infty $, the set ${\bf N}\cup \{ \omega _0 \} $ is countable,
then there exists $k\in {\bf N}\cup \omega _0$ such that the outer
measure of the intersection is positive $\mu ^*(U\cap G_k)>0$, where
$G_k= gr \{ g\in G: ord (g)=k \} $ is the minimal subgroup in $G$
generated by elements of the $k$-th order. Therefore, it is
sufficient to consider all such $G_k$ in the topology inherited from
$G$, $\mu ^*(U\cap G_k)>0$.
\par Mention that if $Y$ is an everywhere dense in $G$ subgroup,
also $card (B\cap T)\ge \sf c$ for some subset $B$ in $G$ for each
$T$ from the base $\Pi $ of neighborhoods of the unit element $e\in
G$, $T\in \Pi $, then $card ((BY)\cap P)\ge \sf c$ for each $P$ open
in $G$.
\par Consider the family $\Xi $ consisting of subgroups $H$ in $G$ and
their automorphisms $s: H\to H$ such that $card (H)\ge {\sf c}$,
$card (s(P)\cap T)\ge {\sf c}$ for each $P$ and $T$ open in $H$ in
the topology inherited from $G$.  Such $H$ exists. For the proof of
their existence take $U$ and  $W$ as above. Let $W_k := \{ g\in W:
ord (g)= k \} $, where $k\in {\bf N}\cup \omega _0$. Then for at
least one $k$ there is satisfied the equality $card (W_k)=card (W)$,
since $card (W)\ge \sf c$. \par Since $card (U)\ge \sf c$ and every
subgroup $\{ g^n: n\in {\bf Z} \} $ generated by a chosen element
$g\in G$ is either finite or countable, then in $W$ there exists a
family $J_{b,k}$, $b\in \Lambda $, $k\in {\bf N}\cup \omega _0$,
$card (\Lambda ) \ge \sf c$, $card (J_b)=card (U)\ge \sf c$. This
family can be chosen such that elements in $J_{b,k}$ would be
algebraically independent: $g\in J_{b,k}$ can not be presented as a
finite product of elements different from it from the set $J_{b,k}$,
moreover, $gr (J_{b,k})\cap gr (J_{d,l})=\{ e \} $ for each $b\ne
d$, since ${\sf q}^{\aleph _0}=\sf q$ for each ${\sf q}\ge \sf c$,
where $J_b=\bigcup_{k\in {\bf N}\cup \omega _0}J_{b,k}$, $\sf q$ and
${\sf c} := card ({\bf R})$ are cardinal numbers, $gr (B)$ denotes
the minimal subgroup in $G$ generated by elements from $B$, every
$g$ from $J_{b,k}$ has the order $ord (g)=k$. Moreover, $J_{b,k}$
can be chosen such that $G\setminus G_{\Lambda }\supset (Y\setminus
\{ e \} )$, where $Y$ is some everywhere dense in $G$ subgroup
$G_{\Lambda } := gr (\bigcup_{b\in \Lambda , k\in {\bf N} \cup
\omega _0} J_{b,k})$, since $G_{\Lambda }$ and $Y$ have
algebraically indeppendent from each other families of generating
elements, $G_{\Lambda }\cap Y= \{ e \} $.
\par Choose every $J_b$ such that $card (J_{b,k_0}\cap P)\ge \sf c$
for each $P$ open in $W$ for at least one $k_0\in {\bf N}\cup \omega
_0$. Make it for every $k_0\in {\bf N}\cup \omega _0$, for which
$card ( \{ g\in W: ord (g)=k_0 \} ) \ge {\sf c} $. Let $\phi :
\Lambda \to \Lambda $ be a bijective mapping from $\Lambda $ onto
$\Lambda $, $\phi (\Lambda (k))=\Lambda (k)$, $\Lambda
=\bigcup_{k\in {\bf N}\cup \omega _0}\Lambda (k)$. Let
$S_{b,k_0}\subset J_{b,k_0}$ and $card (S_{b,k_0})=card
(S_{d,k_0})\ge \sf c$ for each $b, d\in \Lambda (k_0)\subset \Lambda
$. With the help of the transfinite induction take $S_{b,k_0}$ and
bijective mappings $\psi ^d_{b,k_0}: S_{d,k_0}\to S_{b,k_0}$ from
$S_{d,k_0}$ onto $S_{b,k_0}$, put $s(g) := \psi ^d_{\phi
(d),k_0}(g)$ for each $g\in S_{d,k_0}$ and with the help of finite
products of elements in $G$, take $s|_Y=id$, extend $s$ from $Y\cup
\bigcup_{b\in \Lambda (k_0), k_0\in {\bf N}\cup \omega _0}
S_{b,k_0}$ up to the automorphism $s: H\to H$, where $H := gr (Y\cup
\bigcup_{b\in \Lambda (k_0), k_0\in {\bf N}\cup \omega _0} S_b)$
(see also \cite{eng}).

\par The family $\Xi $ order with the help of the relation
$(H_1,s_1)\le (H_2,s_2)$, if $H_1\subset H_2$ and $s_2|_{H_1}=s_1$,
that supplies $\Xi $ with the structure of the directed set. Every
linearly ordered subset ${\cal C} = \{ (H_{\beta }, s_{\beta }):
\beta \in \Lambda ({\cal C}) \} $ in $\Xi $ has an element $(H,s)\in
\Xi $ such that $(H_{\beta },s_{\beta })\le (H,s)$ for each
$(H_{\beta },s_{\beta })\in {\cal C}$, where $H= \bigcup_{\beta \in
\Lambda ({\cal C})} H_{\beta }$, $H\subset G$, $s|_{H_{\beta
}}=s_{\beta }$ for each $\beta \in \Lambda ({\cal C})$, here
$\Lambda ({\cal C})\subset \Lambda $. In view of the Kuratowski-Zorn
lemma in $\Xi $ there exists a maximal element \cite{eng}. Then it
need to be $(G,f)$, since each automorphism $s$ from a subgroup has
an extension up to an automorphism $f$ of the entire group
\cite{neumann,focus}. In view of Lemma 12 the automorphism $f$ is
$({\cal A}_{\mu },{\cal B}(G))$-non-measurable. Since there exists
not less than $2^{\sf c}$ different bijective surjective mappings
$\phi : \Lambda \to \Lambda $, then $card (\Upsilon )\ge 2^{\sf c}$.

\par {\bf 21. Remark.} If $G$ is a Lie group over the field $\bf K$,
$T_eG$ is the Banach space of separable type over the field $\bf K$,
then measures $\mu $ on $G$ and $\nu $ on $T_eG$ with the needed
properties from \S 19 exist (see \cite{dalfom,ludanmat}). If $G$ is
a $C^{\infty }$ over $\bf R$ or $C^{\omega }$ over a local
non-archimedean field Banach-Lie group, then $G$ as the manifold has
the $C^{\infty }$ exponential mapping (see
\cite{bourgralg,bourmnog,kling}).

\section{Non-measurable automorphisms of groups relative to
measures with values in local fields}

\par To avoid misunderstandings we first
remind the basic definitions.
\par {\bf 1. Definitions.} Let $X$ be a completely regular totally
disconnected topological space, let also $\cal R$ be its covering
ring of subsets in $X$, $\bigcup \{ A: A\in {\cal R} \} =X$. We call
the ring separating, if for each two distinct points $x, y \in X$
there exists $A\in \cal R$ such that $x\in A$, $y\notin A$. A
subfamily ${\cal S}\subset \cal R$ is called shrinking, if an
intersection of each two elements from $\cal A$ contains an element
from $\cal A$. If $\cal A$ is a shrinking family, $f: {\cal R}\to
\bf K$, where ${\bf K}=\bf R$ or $\bf K$ is the field with the
non-archimedean norm, then it is written $\lim_{A\in \cal S}
f(A)=0$, if for each $\epsilon >0$ there exists $A_0\in \cal S$ such
that $|f(A)|<\epsilon $ for each $A\in \cal S$ with $A\subset A_0$.
\par A measure $\mu : {\cal R}\to \bf K$ is a mapping with values in
the field $\bf K$ of zero characteristic with the non-archimedean
norm satisfying the following properties:
\par $(i)$ $\mu $ is additive;
\par $(ii)$ for each $A\in \cal R$ the set $\{ \mu (B): B\in {\cal R},
A\subset B \} $ is bounded; \par $(iii)$ if $\cal S$ is the
shrinking family in $\cal R$ and $\bigcap_{A\in \cal S}A =\emptyset
$, then $\lim_{A\in \cal S} \mu (A) = 0$. \par Measures on ${\sf
Bco}(X)$ are called tight measure, where ${\sf Bco}(X)$ is the ring
of clopen (simultaneously open and closed) subsets in $X$.
\par For each $A\in \cal R$ there is defined the norm:
$\| A \| _{\mu } := \sup \{ |\mu (B)|: B\subset A, B\in {\cal R} \}
$. For functions $f: X\to \bf K$ and $\xi : X\to [0,+ \infty )$
define the norm $\| f \|_{\xi }:= \sup \{ |f(x)| \xi (x): x\in X \}
$. Put also $N_{\mu } (x) := \inf \{ \| U \|_{\mu }: x\in U\in {\cal
R} \} $. If a function $f$ is a finite linear combination over the
field $\bf K$ of characteristic functions $\chi _A$ of subsets
$A\subset X$ from $\cal R$, then it is called simple. A function $f:
X\to \bf K$ is called $\mu $-integrable, if there exists a sequence
$f_1, f_2,...$ of simple functions such that there exists
$\lim_{n\to \infty } \| f-f_n \| _{N_{\mu }}=0$.
\par The space $L(\mu )=L(X,{\cal R},\mu ,{\bf K})$ of all $\mu
$-integrable functions is $\bf K$-linear. At the same time
$\int_X\sum_{j=1}^n a_j\chi _{A_j}(x)\mu (dx) := \sum_{j=1}^n a_j\mu
(A_j)$ for simple functions extends onto $L(\mu )$, where $a_j\in
\bf K$, $A_j\in \cal R$ for each $j$.
\par Put ${\cal R}_{\mu } := \{ A: A\subset X, \chi _A \in L(\mu ) \} $.
For $A\in {\cal R}_{\mu }$ let ${\bar \mu }(A) := \int_X\chi
_A(x)\mu (dx)$. \par An automorphism $\phi  $ of a totally
disconnected Hausdorff topological group $G$ is called $\mu
$-non-measurable, if it is $({\cal R}_{\mu }, {\cal
R})$-non-measurable.\par A totally disconnected compact Hausdorff
group $G$ is called $p$-free, if it does not contain any open normal
subgroup of an index divisible by $p$.
\par Let $G$ be a totally disconnected Hausdorff locally compact group,
let also ${\sf B_c}(G)$ be a covering ring of clopen compact subsets
in $G$, suppose that $\mu : {\sf B_c}(G)\to \bf K$ is a
finitely-additive function such that its restriction $\mu |_A$ for
each $A\in {\sf B_c}(G)$ is a tight measure. A measure $\mu : {{\sf
B_c}(G)}\to \bf K$ is called left- (right-)invariant Haar measure,
if $\mu (gA)=\mu (A)$ ($\mu (Ag)=\mu (A)$ respectively) for each
$A\in {\sf B_c}(G)$ and $g\in G$. \par A measure $\eta : {\cal R}\to
\bf K$ is called absolutely continuous relative to a measure $\mu :
{\cal R}\to \bf K$, if there exists a function $f\in L(\mu )$ such
that $\eta (A)=\int_X \chi _A(x) f(x)\mu (dx)$ for each $A\in \cal
R$, denote it by $\eta \preceq \mu $. If $\eta \preceq \mu $ and
$\mu \preceq \eta $, then we say that $\eta $ and $\mu $ are
equivalent $\eta \sim \mu $.

\par {\bf 2. Lemma.} {\it Let $X$ be a Tychonoff (completely regular)
totally disconnected dense in itself topological space with a tight
measure $\mu $ such that for each $x\in X$ the function $N_{\mu
}(x)>0$ is non-negative, moreover, $\mu $ has not atoms in $X$. If
$f: X\to X$ is a bijective epimorphic mapping such that $card
(f(U)\cap V)\ge {\sf c} := card ({\bf R})$ for each open subsets $U$
and $V$ in $X$, then $f$ and $f^{-1}$ are not $({\cal R}_{\mu },{\sf
Bco}(X))$-measurable.}
\par {\bf Proof.} Since $f$ is the bijective mapping from $X$
onto $X$, then $f^{-1}(f(U)\cap V)=U\cap f^{-1}(V)$, consequently,
$card (U\cap f^{-1}(V)\ge \sf c$ for each $U$ and $V$ open in $X$.
In view of Lemma 7.2 \cite{rooij} $ \| \chi _U \|_{N_{\mu }} = \| U
\| _{\mu }$ for each $U\in \cal R$. Due to Lemma 7.5 \cite{rooij}
$N_{\mu }(x)=N_{\bar \mu }(x)$ for each $x\in X$, $L(\mu )=L({\bar
\mu })$ and ${\cal R}_{\bar \mu }={\cal R}_{\mu }$. Theorem 7.6
\cite{rooij} states that $N_{\mu }(x)$ is upper semi-continuous and
for each $\epsilon >0$ the set $\{ x\in A: N_{\mu }(x)\ge \epsilon
\} $ is $R_{\mu }$-compact. If $x_0\in X$, then $0\le N_{\mu }
(x_0)<\infty $, and for each $x_0\in X$ and $r> N_{\mu }(x_0)$ there
exists a neighborhood $P$ of the point $x_0$ such that for each
point $x\in P$ there is accomplished the inequality $N_{\mu }(x)<r$.
Then for each $\epsilon >0$ there exists an open subset $W$ in $X$
such that $A\subset W$ and $N_{\mu }(x) < \epsilon $ for every $x\in
W\setminus A$.
\par For arbitrary clopen subsets $U$ and $V$ of finite
$\mu $ measure take clopen subsets $U_1$ and $U_2$ in $U$, $V_1$ and
$V_2$ in $V$ such that $0< \|U \|_{\mu }<\infty $, $0< \| V \| _{\mu
}<\infty $, $U_1\cap U_2=\emptyset $, $V_1\cap V_2=\emptyset $,
$0<\delta \le \min ( \| U_1 \| _{\mu }, \| U_2 \| _{\mu }, \| V_1 \|
_{\mu }, \| V_2 \| _{\mu }) /s$, $ \| U\setminus (U_1\cup U_2) \|
_{\mu }\le \delta /s$,  $ \| V\setminus (V_1\cup V_2) \| _{\mu }\le
\delta /s$. This is possible, since $N_{\mu }(x)>0$ for each $x\in
X$, while $\mu $ has not atoms.
\par Denote $A:= f^{-1} (U)$, $A_j := f^{-1}(U_j)$. By the supposition of
this Lemma $card (A_j\cap V_k)\ge \sf c$ for each $j, k \in \{ 1, 2
\} $. Suppose that $f$ is the $({\cal R}_{\mu },{\sf Bco}(X))$
measurable mapping. Then there would be $A_j\in L({\mu })$ for $j=1,
2$ and there would be open subsets $B_j$ in $X$ such that
$A_j\subset B_j$ and $N_{\mu }(x) < \delta /s$ for each $x\in
B_j\setminus A_j$ and $j=1, 2$. But $A_1\cap A_2=\emptyset $,
consequently, $\mu (A)=\mu (A_1)+\mu (A_2)$. \par Since $N_{\mu
}(x)>0$ for each $x\in X$, then there can be chosen clopen $U$,
$U_1$ and $U_2$ such that $N_{\mu }|_{(A_j\cap V)}>0$, where
$U_1\cup U_2\subset U$. But $card (A_j\cap P_k)\ge \sf c$ for each
$P_k$ open subset in $V_k$.  Consequently, there exists a countable
sequence of open subsets $Y_n$ in $X$ with $Y_n\supset (A_j\cap
V_k)$ such that $\lim_{n\to \infty } \sup \{ N_{\mu }(x): x\in
Y_n\setminus (A_j\cap V_k) \}  =0$ for given $j, k$,
$Y_n=Y_{n,j,k}$, $j, k\in \{ 1, 2 \} $. Let $C_n =
\bigcap_{s=1}^nY_s$, then $(A_j\cap V_k)\subset C_{n+l}\subset C_n$
for any $n, l\in \bf N$, every $C_n$ is open in $X$,
$C_n=C_{n,j,k}$. \par At the same time $\lim_{n\to \infty } \sup \{
N_{\mu }(x): x\in C_n\setminus (A_j\cap V_k) \} =0$. Since $card
(A_j\cap P)\ge \sf c$ for each subset $P$ open in $X$, then for each
$\epsilon >0$ there exists $n_0\in \bf N$ such that $ Int (\{ x\in
C_n\setminus C_{n+l}: N_{\mu }(x)\ge \epsilon \} ) = \emptyset $ for
each $n>n_0$ and $l\in \bf N$, where $Int (B)$ denotes the interior
of the subset $B$ in $X$. Then $\mu  (A_j\cap V_k) = {\bar \mu
}(cl_X(A_j\cap V_k))= \mu (cl_X(A_j\cap V_k))=\mu (V_k)$, since
$cl_X(A_j\cap V_k)=V_k \in {\cal R}_{\mu }(X)$, where $X$ is the
completely regular space dense in itself, here $cl_X(B)$ denotes the
closure of a subset $B$ in $X$. \par But then it would be
$\sum_{j=1}^2 \mu (A_j\cap V_k) = 2 \mu (V_k)$, that contradicts to
the additivity of the measure: $2\mu (V_k) = \mu (A_1\cap V_k) +\mu
(A_2\cap V_k) = \mu ((A_1\cup A_2)\cap V_k)={\bar \mu }(cl_X(A_1\cup
A_2)\cap V_k))=  \mu (cl_X(A_1\cup A_2)\cap V_k))= \mu (V_k)$, since
the characteristic of the field $\bf K$ is zero, $char ({\bf K})=0$.
Consequently, $f$ is not measurable (see also Theorem 7.12
\cite{rooij}).
\par Applying this proof to $f^{-1}$ instead of $f$ we get that
$f^{-1}$ as well is $\mu $-non-measurable, since $f^{-1}$ satisfies
conditions of the second section of the proof.

\par {\bf 3. Theorem.} {\it Let $G$ be a non-trivial locally
compact Lie group over the non-archimedean local field ${\bf F}$,
${\bf F}\supset {\bf Q_p}$, and $\mu $ be a non-trivial tight Haar
measure on $G$ with values in a local field ${\bf K}\supset {\bf
Q_s}$, where $p$ and $s$ are mutually prime numbers, $(p,s)=1$. Then
the group of its automorphisms $Aut (G)$ has a family of the
cardinality not less than $2^{\sf c}$ of distinct $\mu
$-non-measurable automorphisms on $G$, where ${\sf c}:= card ({\bf
R})$ denotes the cardinality of the continuum.}

\par {\bf Proof.} Since the group $G$ belongs to the class
of smoothness $C^{\omega }$, then for it the Lie algebra $\sf g$
over $\bf F$ is defined. This Lie algebra is the finite-dimensional
space over $\bf F$, $dim_{\bf F}{\sf g}=n\in \bf N$. Then the
additive group for $\sf g$ is $s$-free. In view of the
Monna-Springer theorem 8.4 \cite{rooij} a Haar measure $\nu ^n$ on
it is defined with values in $\bf K$ such that $\nu ^n(B({\bf
F^n},0,1))=1$. It is known, that the Haar measure on ${\sf B_c}(G)$
has not atoms.
\par Take a clopen compact subgroup $W$ in $G$
and an automorphism $\phi $ of the group $G$ from \S 2.13. Without
loss of generality choose $W$ such that for it there is satisfied
the Campbell--Hausdorff formula. In view of Lemma 2 of this section
the automorphism $\phi $ is not $\mu $-measurable.

\par {\bf 4. Corollary.} {\it The family $\Upsilon $ of non-measurable
automorphisms from Theorem 3 has a subfamily $\Omega $ of the
cardinality $card (\Omega ) \ge 2^{\sf c}$ such that every $f\in
\Omega $ being restricted on any one-parameter subgroup over the
field $\bf F$ in $G$ is non-measurable relative to the corresponding
Haar measure on the subgroup with values in $\bf K$.}
\par {\bf Proof.} The exponential mapping
$\exp $ from a neighborhood of zero $V_0$ in the algebra $\sf g$
onto a neighborhood $U_e$ of the unit element in $G$ induces the
image of the measure $\nu ^n_{\exp }$ on $U_e$, where $n$ is the
dimension of $\sf g$ as the linear space over the field $\bf F$,
$\nu ^n$ is the Haar measure on $\bf F^n$ as the additive group,
moreover, it has not atoms. \par Theorem 7.34 in \cite{rooij} states
that if there are two measures $\lambda $ and $\zeta $ on a covering
ring $\cal R$ of a topological completely regular totally
disconnected space $X$, then the following two conditions are
equivalent: $(\alpha )$ there exists a locally $\zeta $-integrable
function $h$ such that $\lambda (dx)= h(x)\zeta (dx)$; $(\beta )$
for each $x\in X$ there exists $b\in \bf K$ with $N_{\lambda -b\zeta
}(x)=0$. Consequently, the measure $\nu ^n_{\exp }$ is equivalent to
the restriction of the Haar measure $\mu $ on $U_e$, since $\exp $
is the locally bijective of class $C^{\omega }$ mapping. Then the
Haar measure $\nu $ on $\bf F$ also induces the measure $\eta _g$ on
the one-parameter subgroup $g_W = \{ g^t: t\in {\bf F}, |t|<\epsilon
\} $, $0<\epsilon $, where $g\in g_W$. This measure $\eta _g$ is
equivalent to the Haar measure $\mu _g$ on $g_W$.
\par Then from \S \S 2.13, 3.3 and Lemma 3.2 the statement of this
corollary follows.

\par {\bf 5. Theorem.} {\it Let $\sf g$ be a non-trivial
Lie algebra finite-dimensional over the field $\bf F$ with a measure
$\mu $ equal to the non-trivial $\bf K$-valued Haar measure on an
additive group of $\sf g$. Then the algebra $\sf g$ has the family
of the cardinality $2^{\bf c}$ of $\mu $-non-measurable
automorphisms.}
\par {\bf Proof.} Take any algebraic automorphism
$\phi $ of the field $\bf F$ from the proof of Theorem 2.13. Since
$\sf g$ is finite-dimensional over the field $\bf F$, then the
non-trivial Haar measure $\mu $ on $\sf g$ as the additive group for
$\sf g$ with values in $\bf K$ is equivalent to the measure $\nu ^n$
(see Definitions 1). The automorphism $\phi $ extends up to an
automorphism of the algebra: $\phi (a_jv_j)=\phi (a_j)v_j$, $\phi
(a_1v_1+....+a_mv_m)= \phi (a_1)v_1+...+\phi (a_m)v_m$, $\phi
([a_kv_k,a_jv_j]) = [\phi (a_k)v_k, \phi (a_j)v_j]$ for each $a_j\in
\bf F$, $k, j=1,...,m$, where $v_1,...,v_m$ is the basis of
generators in $\sf g$  (see Definitions 2.1). Therefore, due to
Lemma 2 the automorphism $\phi $ of the algebra $\sf  g$ is $({\sf
B_c}_{\mu },{\sf B_c }({\sf g}))$-non-measurable. The family of such
different automorphisms of the algebra $\sf g$ has the cardinality
$2^{\sf c}$.

\par {\bf 6.} Let $G$ be a $C^{\omega }$ non-trivial Lie group
over a non-archimedean local field $\bf F$, moreover, $G$ be
complete as an uniform space. Suppose that $\mu : {\cal R}(G)\to \bf
K$ is finitely-additive and there exists a clopen subgroup $W$ in
$G$ such that the restriction $\mu |_W$ is the tight non-trivial
measure on ${\sf Bco}(W)$ with values in $\bf K$ such that $N_{\mu
}(g)>0$ for each $g\in G$, moreover, $\mu $ has not any atoms, where
${\sf Bco}(G)\supset {\cal R}(G)\supset {\sf Bco}(W)$, ${\cal R}(G)$
is a covering ring for $G$. Let $\exp : V\to W$ be the exponential
mapping for $G$ as the analytic $C^{\omega }$ manifold over $\bf F$
from an open neighborhood $V$ of zero in $T_eG$ on $W$. Also suppose
that a tight $\bf K$-valued measure $\nu $ on ${\cal R}(T_eG)$ is
such that $\nu (J)=\mu (\exp (J))$ for each $J\in {\sf Bco}(V)$,
where $T_eG$ is the linear space of separable type over the field
$\bf F$, ${\sf Bco}({\sf g})\supset {\cal R}({\sf g})\supset {\sf
Bco}(V)$, ${\cal R}({\sf g})$ is the covering ring for $\sf g$. Let
as in Theorem 3 $p$ and $s$ be two mutually prime numbers,
$(p,s)=1$.
\par Suppose that $G$ has a clopen subgroup $W$,
$\exp : V\to W$, moreover, in $W$ there exists an everywhere dense
subgroup $S$ such that the restriction $\ln |_S$ corresponds to the
Campbell-Hausdorff formula, where $\ln $ is the inverse mapping to
$\exp $.
\par Let $\pi _v: T_eG\to {\bf F}v$ be a $\bf F$-linear projection
operator, where $\nu _v$ is equivalent to the Haar measure on $\bf
F$, where $v$ is the non-zero vector of the tangent space $v\in
T_eG$, $\nu _v(J)=\nu (\pi _v^{-1}(J))$ for each $J\in {\sf
B_c}({\bf F})$.

\par {\bf Theorem.} {\it Then such group $G$ has a family
$\Upsilon $ of $({\cal R}(G)_{\mu },{\cal R}(G))$-non-measurable
automorphisms of the cardinality not less than $2^{\sf c}$, $card
(\Upsilon ) \ge 2^{\sf c}$. Moreover, $\Upsilon $ has a subfamily
$\cal P$ of automorphisms $f$, restrictions of which on
one-parameter over $\bf F$ local subgroups $\{ \exp (xv):
|x|<\epsilon \} $ in $S$ are non-measurable relative to the
corresponding $\bf K$-valued Haar measure on $\{ \exp (xv):
|x|<\epsilon \} $.}
\par {\bf Proof.} The image $\nu $ on $V$ of the measure $\mu $ with the
help of the mapping $\ln $, $\nu (B)= \mu (\exp (B))$ for every
$B\in {\sf Bco}(V)$ is extendable up to the tight measure on the
corresponding covering ring ${\cal R}({\sf g})$, $\nu (B) :=
\sum_{j=1}^{\infty } \nu ((B-h_j)\cap V)s^j$ for every $B\in {\cal
R}({\sf g})$, where $\{ (V+h_j): j\in {\bf N}, h_j\in {\sf g} \} $
is the covering for $\sf g$, since $\sf g$ by the condition has the
separable type over $\bf F$, while the field $\bf F$ is separable
and locally compact. As ${\cal R}({\sf g})$ we can take the minimal
ring generated by $\bigcup_{j=1}^{\infty } {\sf Bco}(V+h_j)$.
Therefore, $\nu $ has tight measures as the projections $\nu
_{\omega ({\bf F^n})}$ on $\omega ({\bf F^n})$ for each embedding
$\omega : {\bf F^n} \hookrightarrow \sf g$ as the $\bf F$-linear
space, $\nu _{\omega ({\bf F^n})}(B) = \nu (\pi ^{-1}(B))$ for each
$B\in {\sf B_c}(\omega ({\bf F^n}))$, where $\pi : {\sf g}\to \omega
({\bf F^n})$ is the projection operator.
\par The operator $\pi $ is $\bf F$-linear and it exists, since
$\omega ({\bf F^n})$ is finite-dimensional over $\bf F$, while the
field $\bf F$ is locally compact (see Theorems 5.13 and 5.16
\cite{rooij} and \cite{nari}). The image $\nu _{\omega ({\bf
F})}|_{V\cap \omega ({\bf F})}$ with the help of $\exp $ generates
the measure on the local one-parameter subgroup in $S$, $\mu _g(B)=
\nu _{\omega ({\bf K})}(\ln (B))$ for every $B\in {\cal R}(g_W\cap
W)$, which extends up to a tight measure $\mu _g$ on $g_W$.
\par Take an automorphism $f$ of the group $G$ from \S 2.19. In view of
Lemma 3.2 and $N_{\mu }(g)>0$ for each $g\in G$ we get that $f$ is
$\mu $-non-measurable. The non-measurability of restrictions of $f$
on one-parameter subgroups relative to the  Haar measures on them
follows from the properties of $\phi $ again due to Lemma 3.2 and \S
2.13. The families $\Upsilon $ and $\cal P$ of such different
automorphisms of the group $G$ due to \S 2.19 have the cardinalities
not less than $2^{\sf c}$.

\par {\bf 7. Theorem.} {\it Let $G$ be an infinite topological
totally disconnected dense in itself Hausdorff group with a
non-trivial tight measure $\mu $ on $G$ having no any atom,
moreover, $N_{\mu }(g)>0$ for each $g\in G$, while $G$ is complete
as the uniform space, $card (U)\ge \sf c$ for each open $U$ in $G$.
Then there exists a family $\Upsilon $ of $\mu $-non-measurable
distinct automorphisms of the group $G$ of the cardinality $card
(\Upsilon )\ge 2^{\sf c}$.}
\par {\bf Proof.} Take an automorphism $\phi $
of the group $G$ from \S 2.20. In view of Lemma 3.2 it is $\mu
$-non-measurable. The family of such distinct automorphisms of the
group $G$ due to \S 2.20 has the cardinality not less than $card
(\Upsilon )\ge 2^{\sf c}$.

\par {\bf 8. Remark.} If $G$ is a Lie group over the field $\bf F$
of the class $C^{\omega }$, its tangent space $T_eG$ is a Banach
space of separable type over the field $\bf F$, then there exist
measures $\mu $ on $G$ and $\nu $ on $T_eG$ with the desired
properties from \S 3.6 (see \cite{lujmsqim}).

\end{document}